\documentclass{article}
\usepackage{amssymb}
\usepackage{latexsym}
\usepackage{amsmath}
\usepackage{mathrsfs}
\usepackage{color}
\usepackage{CJK}
\usepackage{graphicx}

\newtheorem{Thm}{Theorem}[section]

\newtheorem{Lem}{Lemma}[section]
\newtheorem{Prop}{Proposition}[section]
\newtheorem{Rmk}{Remark}[section]
\newtheorem{Def}{Definition}[section]

\numberwithin{equation}{section}
\newenvironment{prooff}{\medskip\par\noindent{\bf Proof.}\ }{\qquad
\raisebox{-0.5mm}{\rule{1.5mm}{4mm}}\vspace{6pt}}

\newcommand{\bbrn}{\mathbb{R}^N}
\newcommand{\bbr}{\mathbb{R}}
\newcommand{\h}{H^1_0(\Omega)}
\newcommand{\lp}{L^p(\bbrn)}
\newcommand{\bbn}{\mathbb{N}}

\newcommand{\ve}{\varepsilon}

\begin{document}

\title
{\Large\bf On the Br\'ezis--Nirenberg problem with nonhomogeneous Dirichlet boundary conditions}

\author{Yuanze Wu$^a$\thanks{Corresponding
author. E-mail address: wuyz850306@cumt.edu.cn (Yuanze Wu)}\quad
Tsung-fang Wu$^b$\thanks{E-mail address:  tfwu@nuk.edu.tw (Tsung-fang Wu)}\quad
Zeng Liu$^c$\thanks{E-mail address: luckliu@163.com (Zeng Liu)}\\
\footnotesize$^a${\em  College of Sciences, China University of Mining And Technology, Xuzhou 221116 P.R. China}\\
\footnotesize$^b${\em  Department of Applied Mathematics, National University of Kaohsiung, Kaohsiung 811, Taiwan }\\
\footnotesize$^c$ {\em Department of Mathematics, Soochow
University, Suzhou 215006, P.R. China}}

\date{}
\maketitle

\date{} \maketitle

\noindent{\bf Abstract:} In this paper, we study the decomposition of Nehari manifold for the Br\'ezis--Nirenberg problem with nonhomogeneous Dirichlet boundary conditions.  By using this result, the Lusternik--Schnirelman category and the minimax principle, we establish a multiple result (four solutions) for the Br\'ezis--Nirenberg problem with nonhomogeneous Dirichlet boundary conditions.
\vspace{6mm}


\vspace{6mm}\noindent {\bf AMS} Subject Classification 2010: 35B38; 35J20; 35J25.

\section{Introduction}
In this paper, we consider the following elliptic problem:
$$
\left\{\aligned-\Delta u&=\lambda u+u^{2^*-1}&\text{ in }\Omega,\\
u&>0&\text{ in }\Omega,\\
u&=\mu g&\text{ on }\partial\Omega,\endaligned\right.\eqno{(\mathcal{P}_{\lambda,\mu})}
$$
where $\Omega\subset\bbrn(N\geq3)$ is a bounded domain with smooth boundary, $2^*:=2N/(N-2)$ is the critical Sobolev exponent, $\lambda, \mu\geq0$ are two parameters, $g(x)\in H^{1/2}(\partial\Omega)\cap C(\partial\Omega)$.

By the classical regularity theory (cf. Gilbarg and Trudinger \cite{GT98}),  each solution of
Problem~$(\mathcal{P}_{\lambda,\mu})$ is $C^1$, which implies that $(\mathcal{P}_{\lambda,\mu})$ has no solution if $\min\{g(x), 0\}\not=0$ for $\mu>0$.  On the other hand, it is easy to see that $(\mathcal{P}_{\lambda,\mu})$ with $\mu>0$ and $g\equiv0$ is equivalent to $(\mathcal{P}_{\lambda,0})$.  So, it is natural to discuss $(\mathcal{P}_{\lambda,\mu})$ under the following assumption on $g(x)$
\begin{itemize}
\item[$(G)$] $g(x)\geq0$ and $g(x)\not\equiv0$.
\end{itemize}

$(\mathcal{P}_{\lambda,0})$ is called as the Br\'ezis--Nirenberg problem, since it has been first studied by Br\'ezis and Nirenberg in their nice paper \cite{BN83}.  In that paper, they have proved that $(\mathcal{P}_{\lambda,0})$ has a solution for $N\geq4$ if and only if $\lambda\in(0, \lambda_1)$, where $\lambda_1$ is the first eigenvalue of $(-\Delta, \h)$.  They have also proved that there exists $\lambda^*\in(0, \lambda_1)$ such that $(\mathcal{P}_{\lambda,0})$ has a solution for $N=3$ when $\lambda\in(\lambda^*, \lambda_1)$ and no solution for $N=3$ when $\lambda\geq\lambda_1$.  Moreover, when $\Omega$ is a ball in $\bbr^3$, $\lambda^*=\frac{\lambda_1}{4}$ and $(\mathcal{P}_{\lambda,0})$ had no solution for $\lambda\in(0, \lambda^*]$ (see also in \cite{BN83}).  Since then, many papers have been devoted to the Br\'ezis--Nirenberg problem, see for example Arioli et al. \cite{AGGS08}, Clapp and Weth \cite{CW05}, Chen et al. \cite{CSZ12}, Schechter and Zou \cite{SZ10}, Zhang \cite{Z89} and the references therein.

When $\mu>0$ and $(G)$ hold, $(\mathcal{P}_{\lambda,\mu})$ is a kind of the so-called elliptic problems with nonhomogeneous Dirichlet boundary conditions.  Since $(\mathcal{P}_{\lambda,\mu})$ links closed to $(\mathcal{P}_{\lambda,0})$, the Br\'ezis--Nirenberg problem, we call $(\mathcal{P}_{\lambda,\mu})$ as the Br\'ezis--Nirenberg problem with nonhomogeneous Dirichlet boundary conditions.

The elliptic problems with nonhomogeneous Dirichlet boundary conditions have been studied widely in the past decade.  To our best knowledge, many papers have been devoted to the existence of infinitely many nontrivial solutions for subcritical cases, see for example Salvatore \cite{S01}, Bolle et al. \cite{BGT00}, Candela et al. \cite{CSS02}, Hu \cite{H04} and the references therein, and there are few results about the critical case, for example $(\mathcal{P}_{\lambda,\mu})$.  Thus, the purpose of this paper is to study the existence and multiplicity of solutions for the elliptic problems with nonhomogeneous Dirichlet boundary conditions and critical Sobolev exponent.  More precisely, the solutions for $(\mathcal{P}_{\lambda,\mu})$.

It is well known that if $v$ is a solution of the following elliptic problem
$$
\left\{\aligned-\Delta v&=\lambda(v+\mu\varphi)+(v+\mu\varphi)^{2^*-1}&\text{ in }\Omega,\\
v&>0&\text{ in }\Omega,\\
v&=0&\text{ on }\partial\Omega,\endaligned\right.\eqno{(\mathcal{Q}_{\lambda,\mu})}
$$
then $u=v+\varphi$ is a solution of $(\mathcal{P}_{\lambda,\mu})$, where $\varphi$ satisfies
$$
\left\{\aligned\Delta\varphi&=0&\text{ in }\Omega,\\
\varphi&=g&\text{ on }\partial\Omega,\endaligned\right.
$$
and vice versa by the maximum principle.  Therefore, the study of the existence and multiplicity of solutions for $(\mathcal{P}_{\lambda,\mu})$ is equivalent to the study of the existence and multiplicity of solutions for $(\mathcal{Q}_{\lambda,\mu})$.

Clearly, $(\mathcal{Q}_{\lambda,\mu})$ has a variational structure in $\h$ and the corresponding functional of $(\mathcal{Q}_{\lambda,\mu})$, defined on $\h$, is given by
$$
I_{\lambda,\mu}(v)=\frac12\int_\Omega|\nabla v|^2dx-\frac\lambda2\int_\Omega|v+\mu\varphi|^2dx-\frac1{2^*}\int_\Omega|v+\mu\varphi|^{2^*}dx.
$$
Hence, we can use the variational method to study the existence and multiplicity of solutions for $(\mathcal{Q}_{\lambda,\mu})$.

The Nehari manifold is a useful tool in proving the existence and multiplicity of solutions for elliptic problems with homogeneous Dirichlet boundary conditions and variational structure (cf. Sun and Li \cite{SL08} and Wu \cite{W10}).  Let $\mathcal{N}_{\lambda,\mu}:=\{u\in \h\backslash\{0\}: I_{\lambda,\mu}'(u)u=0\}$ be the Nehari manifold of $I_{\lambda,\mu}$.  It is well known that the Nehari manifold $\mathcal{N}_{\lambda,\mu}$ is closely linked to the behavior of the fibering maps, which is given by  $T_{\lambda,\mu,u}(t):=I_{\lambda,\mu}(tu)$, $t>0$.  The fibering map has been introduced by Dr\'abek and Pohozaev in \cite{DP97} and has also been studied by Brown and Wu \cite{BW09}, Brown and Zhang \cite{BZ03}.  Moreover,  $T_{\lambda,\mu,u}'(t)=0$ is equivalent to $tu\in\mathcal{N}_{\lambda,\mu}$.  In particular, $T_{\lambda,\mu,u}'(1)=0$ if and only if $u\in\mathcal{N}_{\lambda,\mu}$.  Since $T_{\lambda,\mu,u}\in C^2(\bbr^+,\bbr)$ for every $0\not=u\in\h$, it is natural to split the Nehari manifold $\mathcal{N}_{\lambda,\mu}$ into three parts:
$$
\aligned
\mathcal{N}_{\lambda,\mu}^+:=&\{u\in\mathcal{N}_{\lambda,\mu} : T_{\lambda,\mu,u}''(1)>0\};\\
\mathcal{N}_{\lambda,\mu}^-:=&\{u\in\mathcal{N}_{\lambda,\mu} : T_{\lambda,\mu,u}''(1)<0\};\\
\mathcal{N}_{\lambda,\mu}^0:=&\{u\in\mathcal{N}_{\lambda,\mu} : T_{\lambda,\mu,u}''(1)=0\}.
\endaligned
$$
The crucial point of using the method of the fibering maps (or the tool of the Nehari manifold) is to show that the unit sphere of $\h$ is homeomorphous to $\mathcal{N}_{\lambda,\mu}^-$.  However, when $\mu>0$, $|u+\mu\varphi|^{2^*}$ can not be controlled by some homogeneous terms of $|u|$. This brings about the usual arguments of proving the homeomorphism between the unit sphere of $\h$ and $\mathcal{N}_{\lambda,\mu}^-$ are invalid in dealing with $(\mathcal{Q}_{\lambda,\mu})$ for $\mu>0$.  Another difference in using the method of the fibering maps to find solutions of $(\mathcal{Q}_{\lambda,\mu})$ with $\mu>0$ is that even though $v_{\lambda,\mu}$ is a minimizer of $I_{\lambda,\mu}(v)$ on $\mathcal{N}_{\lambda,\mu}^\pm$, it is not trivial to know $|v_{\lambda,\mu}|$ is also a minimizer of $I_{\lambda,\mu}(v)$ on $\mathcal{N}_{\lambda,\mu}^\pm$.  In order to prove the multiplicity of solutions (more precisely, four solutions) for $(\mathcal{Q}_{\lambda,\mu})$ by using the method of fibering maps, we need to overcome these difficulties, which require us to make some careful analysis and complicated estimates.

The first result obtained in this paper is
\begin{Thm}\label{thm003}
Assume $(G)$ holds and $N\geq4$.  Then we have the following.
\begin{itemize}
\item[$(a)$] $(\mathcal{P}_{\lambda,\mu})$ has no solution for $\lambda\in[\lambda_1, +\infty)$ and $\mu\geq0$.
\item[$(b)$] For every $\lambda\in(0, \lambda_1)$, there exist $0<\mu_\lambda^{*}\leq\mu^*(\lambda)<+\infty$ such that $(\mathcal{P}_{\lambda,\mu})$ has two solutions for $\mu\in(0, \mu_\lambda^{*})$, one solution for $\mu\in[\mu_\lambda^{*}, \mu^*(\lambda)]$ and no solutions for $\mu>\mu^*(\lambda)$.
\end{itemize}
\end{Thm}

It is well know that the number of solutions for elliptic problems on bounded domains are effected by the topology of the domain, see for example He and Yang \cite{HY09} and Wu \cite{W10}.  Inspired by these papers, we obtain the following result.
\begin{Thm}\label{thm004}
Assume $(G)$ holds and $N\geq4$.  If the domain $\Omega$ satisfies
\begin{itemize}
\item[$(D)$] There exists $\delta_0\in(0, 1)$ such that $\overline{B}_{\frac1\delta_0}\backslash B_{\delta_0}\subset \Omega$ and $B_{\frac{\delta_0}2}\cap \Omega=\emptyset$, where $B_r:=\{x\in\bbrn:|x|<r\}$.
\end{itemize}
Then there exist $\lambda_*\in(0, \lambda_1)$ and $\mu_{*}>0$ such that
$(\mathcal{P}_{\lambda,\mu})$ has at least three solutions for $\lambda\in(0, \lambda_*)$ and $\mu\in(0, \mu_{*})$.  Furthermore, if $\delta_0$ is small enough, then
$(\mathcal{P}_{\lambda,\mu})$ has at least four solutions for $\lambda\in(0, \lambda_*)$ and $\mu\in(0, \mu_{*})$.
\end{Thm}

In this paper, we will always denote the usual norm in $\h$ and $\lp$ by $\|\cdot\|$ and $\|\cdot\|_p$, respectively.  $C$ will be indiscriminately used to denote various positive constants.

\section{Nonexistence results for $(\mathcal{Q}_{\lambda,\mu})$}
In this section, we will discuss the nonexistence results for $(\mathcal{Q}_{\lambda,\mu})$.  For $\lambda\in(0, \lambda_1)$, we consider an operator which is defined on $\h\times\bbr$ as follows
$$
\mathcal{H}_\lambda(v,\mu)=-\Delta v-\lambda(v+\mu\varphi)-(v+\mu\varphi)^{2^*-1}.
$$
It is clearly that $\mathcal{H}_{\lambda}(0,0)=0$ and $\frac{d\mathcal{H}_{\lambda}}{dv}(0,0)=id-\lambda\mathcal{F}$, where $\mathcal{F}$ is the imbedding map from $\h$ to $L^2(\Omega)$.  Since $\lambda\in(0, \lambda_1)$, by the implicit
function theorem, there exist $\mu^{0}(\lambda)>0$ and a continuous curve $(v_\mu,\mu)$ in $\h\times\bbr$ with $v_0=0$ such that $\mathcal{H}_{\lambda}(v_\mu,\mu)=0$ for all
$0\leq\mu\leq\mu^{0}(\lambda)$.  For every $\lambda\in(0, \lambda_1)$, let
$$
\mu^*(\lambda):=\sup\{\mu: \mu>0\text{ and }(\mathcal{Q}_{\lambda,\mu})\text{ has a solution}\}.
$$
Then it is easy to see that $\mu^*(\lambda)\geq\mu^{0}(\lambda)>0$.  The next lemma implies $\mu^*(\lambda)<+\infty$ for every $\lambda\in(0, \lambda_1)$.
\begin{Lem}\label{lem005}
Assume $(G)$ holds.  Then $\mu^*(\lambda)<+\infty$ for every $\lambda\in(0, \lambda_1)$.
\end{Lem}
\begin{prooff}
Assuming the contrary.  Then there exists $\lambda_0\in(0, \lambda_1)$ such that $\mu^*(\lambda_0)=+\infty$.  This means that there exists $\{\mu_n\}$ satisfying $\mu_n\to+\infty$ as $n\to+\infty$ such that $(\mathcal{Q}_{\lambda_0,\mu_n})$ has a solution $u_{\mu_n}$ for every $n\in\bbn$.  Since $u_{\mu_n}$ is a solution of $(\mathcal{Q}_{\lambda_0,\mu_n})$, multiplying $e_1$, the eigenfunction of $\lambda_1$, on $(\mathcal{Q}_{\lambda_0,\mu_n})$ with both sides and integrate, we obtain
$$
(\lambda_1-\lambda_0)\int_\Omega u_{\mu_n} e_1dx\geq\lambda_0\mu_n\int_\Omega\varphi e_1dx+\int_{\mathcal{A}_{\lambda_0,n}}u_{\mu_n}^{2^*-1}e_1dx,
$$
where $\mathcal{A}_{\lambda_0,n}:=\{x\in\Omega: u_{\mu_n}>(\lambda_1-\lambda_0)^{\frac{1}{2^*-2}}\}$.
On the other hand,
$$
(\lambda_1-\lambda_0)\int_\Omega u_{\mu_n} e_1dx\leq(\lambda_1-\lambda_0)\int_{\mathcal{A}_{\lambda_0,n}} u_{\mu_n} e_1dx+(\lambda_1-\lambda_0)^{\frac{2^*-1}{2^*-2}}\int_\Omega e_1dx.
$$
Therefore,
\begin{equation}\label{eq016}
\int_{\mathcal{A}_{\lambda_0,n}}\big((\lambda_1-\lambda_0)u_{\mu_n}-u_{\mu_n}^{2^*-1}\big)e_1dx+\int_\Omega
(\lambda_1-\lambda_0)^{\frac{2^*-1}{2^*-2}}e_1dx-\lambda_0\mu_n\int_\Omega\varphi e_1dx\geq0.
\end{equation}
Note that $\int_{\mathcal{A}_{\lambda_0,n}}\big((\lambda_1-\lambda)u_{\mu_n}-u_{\mu_n}^{2^*-1}\big)e_1dx<0$.  Thus \eqref{eq016} is impossible for $n$ large enough, since $\int_\Omega\varphi e_1dx>0$ by $(G)$ and the maximum principle.
\end{prooff}

We close this section by the following proposition.
\begin{Prop}\label{prop01}
Assume $(G)$ holds.  Then
\begin{itemize}
\item[$(1)$] $(\mathcal{Q}_{\lambda,\mu})$ has no solution for $\lambda\in[\lambda_1, +\infty)$ and $\mu\geq0$.
\item[$(2)$] For every $\lambda\in(0, \lambda_1)$,  $(\mathcal{Q}_{\lambda,\mu})$ has one solution for $\mu\in(0, \mu^*(\lambda)]$ and no solutions for $\mu>\mu^*(\lambda)$.
\end{itemize}
\end{Prop}
\begin{prooff}
Multiplying $e_1$ on $(\mathcal{Q}_{\lambda,\mu})$ with both sides and integrate, it is easy to see that the conclusion $(1)$ holds since $\varphi\geq0$ under $(G)$ by the maximum principle.  For every
$\lambda\in(0, \lambda_1)$,  by the definition of $\mu^*(\lambda)$, we can use the method of sub- and supper-solutions (cf. Ca\~nada et al. \cite{CDG97}) in a standard way to show that $(\mathcal{Q}_{\lambda,\mu})$ has a solution for every $\mu\in(0, \mu^*(\lambda))$.  By taking a limit, $(\mathcal{Q}_{\lambda,\mu})$ also has a solution for $\mu=\mu(\lambda)^*$.  Thus, conclusion $(2)$ follows immediately from Lemma~\ref{lem005}.
\end{prooff}

\section{Properties of the Nehari manifold}
In this section, we will study the decomposition and properties of the Nehari manifold $\mathcal{N}_{\lambda,\mu}$ and the fibering map $T_{\lambda,\mu,v}(t)$.  We start by
\begin{Lem}\label{lem001}
For every $0<\lambda<\lambda_1$, there exists $\mu_\lambda^0>0$ such that, for every $\mu\in(0, \mu_\lambda^0)$, we have
\begin{itemize}
\item[$(1)$] If $v\not=0$ and $\lambda\mu\int_\Omega\varphi vdx+\mu^{2^*-1}\int_\Omega\varphi^{2^*-1}vdx\leq0$ then there exists unique $t_{\lambda,\mu}^-(v)>0$ such that $t_{\lambda,\mu}^-(v)v\in\mathcal{N}_{\lambda,\mu}^-$.
\item[$(2)$] If $v\not=0$ and $\lambda\mu\int_\Omega\varphi vdx+\mu^{2^*-1}\int_\Omega\varphi^{2^*-1}vdx>0$ then there exists unique $t_{\lambda,\mu}^-(v)>t_{\lambda,\mu}^+(v)>0$ such that $t_{\lambda,\mu}^+(v)v\in\mathcal{N}_{\lambda,\mu}^+$ and $t_{\lambda,\mu}^-(v)v\in\mathcal{N}_{\lambda,\mu}^-$. Moreover, $I_{\lambda,\mu}(t_{\lambda,\mu}^+(v)v)=\min_{0\leq t\leq t_{\lambda,\mu}^-(v)}I_{\lambda,\mu}(tv)$.
\item[$(3)$] $\mathcal{N}_{\lambda,\mu}^{0}=\emptyset$.
\item[$(4)$] $\text{dist}(\mathcal{N}_{\lambda,\mu}^-,\mathcal{N}_{\lambda,\mu}^+):=\inf\{\|v-z\|:v\in\mathcal{N}_{\lambda,\mu}^-, z\in\mathcal{N}_{\lambda,\mu}^+\}>C(1-\frac{\lambda}{\lambda_1})^{\frac{1}{2^*-2}}>0$.
\end{itemize}
\end{Lem}
\begin{prooff}
Clearly, for every $0\not=v\in\h$,
$$
T_{\lambda,\mu,v}'(t)=t\|v\|^2-\lambda\int_\Omega(tv+\mu\varphi)vdx-\int_\Omega|tv+\mu\varphi|^{2^*-2}(tv+\mu\varphi)vdx
$$
and
$$
T_{\lambda,\mu,v}''(t)=\|v\|^2-\lambda\|v\|_2^2-(2^*-1)\int_\Omega|tv+\mu\varphi|^{2^*-2}v^2dx.
$$
Since $\lambda<\lambda_1$, there exists $\mu_\lambda^1>0$ such that $(1-\frac{\lambda}{\lambda_1})S-(2^*-1)2^{2^*-2}\mu^{2^*-2}\|\varphi\|_{2^*}^{2^*}>0$ for $\mu\in(0, \mu_\lambda^1)$, where $S$ is the best Sobolev embedding constant given by $S=\inf_{u\in\h\backslash\{0\}}\frac{\|u\|^2}{\|u\|_{2^*}^2}$.  Let
$$
t_0(v)=(\frac{\|v\|^2-\lambda\|v\|_2^2-(2^*-1)2^{2^*-2}\mu^{2^*-2}\int_\Omega|\varphi|^{2^*-2}v^2dx}
{(2^*-1)2^{2^*-2}\|v\|_{2^*}^{2^*}})^{\frac{1}{2^*-2}},
$$
then $T_{\lambda,\mu,v}''(t)>0$ for $t\in(0, t_0(v))$.  For the sake of clarity, we divide the following proof into several claims.

{\bf Claim~1:}\quad For every $\lambda\in(0, \lambda_1)$, there exists $\mu_\lambda^2\in(0, \mu_\lambda^1)$ such that $T_{\lambda,\mu,v}'(t_0(v))>0$ for $\mu\in(0, \mu_\lambda^2)$.

Indeed, if $\lambda\mu\int_\Omega\varphi vdx+\mu^{2^*-1}\int_\Omega\varphi^{2^*-1}vdx\leq0$, then $T_{\lambda,\mu,v}'(0)\geq0$.  Since $T_{\lambda,\mu,v}''(t)>0$ for $\lambda\in(0, \lambda_1)$, $\mu\in(0, \mu_\lambda^1)$ and $t\in(0, t_0(v))$, $T_{\lambda,\mu,v}'(t_0(v))>0$ for $\lambda\in(0, \lambda_1)$ and $\mu\in(0, \mu_\lambda^1)$.  On the other hand, if $\lambda\mu\int_\Omega\varphi vdx+\mu^{2^*-1}\int_\Omega\varphi^{2^*-1}vdx>0$, then $T_{\lambda,\mu,v}'(0)<0$.  In this case, by directly calculation, we have
$$
\aligned
T_{\lambda,\mu,v}'(t_0(v))&=t_0(v)(\|v\|^2-\lambda\|v\|_2^2)-\lambda\mu\int_\Omega\varphi vdx-\int_\Omega|t_0(v)v+\mu\varphi|^{2^*-2}(t_0(v)v+\mu\varphi)vdx\\
\geq&t_0(v)(\|v\|^2-\lambda\|v\|_2^2)-\lambda\mu\int_\Omega\varphi vdx-2^{2^*-2}(t_0(v)^{2^*-1}\|v\|_{2^*}^{2^*}+\int_\Omega|v||\mu\varphi|^{2^*-1}dx)\\
\geq&t_0(v)(1-\frac{1}{2^*-1})(\|v\|^2-\lambda\|v\|_2^2)-\lambda\mu\int_\Omega\varphi vdx-2^{2^*-2}\int_\Omega|v||\mu\varphi|^{2^*-1}dx\\
\geq&(C(1-\frac{\lambda}{\lambda_1}-C_1\mu^{2^*-2})^{\frac{1}{2^*-2}}(1-\frac{\lambda}{\lambda_1})-C_1\mu^{2^*-1}-C_1\lambda\mu)\|v\|_{2^*}.
\endaligned
$$
Therefore, for every $\lambda\in(0, \lambda_1)$, there exists $\mu_\lambda^2\in(0, \mu_\lambda^1)$ such that $T_{\lambda,\mu,v}'(t_0(v))>0$ for $\mu\in(0, \mu_\lambda^2)$.

{\bf Claim~2:}\quad For every $\lambda\in(0, \lambda_1)$, there exists $\mu_\lambda^3\in(0, \mu_\lambda^2)$ such that if $t\geq t_0(v)$, $\mu\in(0, \mu_\lambda^3)$ and $T_{\lambda,\mu,v}'(t)=0$ then $T_{\lambda,\mu,v}''(t)<0$.

In fact, if not, then there exist $\lambda_0$ and a sequence $\mu_n\to0$ as $n\to\infty$ such that $T_{\lambda_0,\mu_n,v}'(t_n(v))=0$ and $T_{\lambda_0,\mu_n,v}''(t_n(v))\geq0$ for some $t_n(v)\geq t_0(v)$.
This, together with the Young inequality, implies
\begin{equation}\label{eq999}
\aligned
0\geq&(2^*-1)t_n(v)^2\int_\Omega|t_n(v) v+\mu_n\varphi|^{2^*-2}v^2dx-\lambda_0\mu_nt_n(v)\int_\Omega\varphi vdx\\
&-\int_\Omega|t_n(v) v+\mu_n\varphi|^{2^*-2}(t_n(v)^2 v+t_n(v)\mu_n\varphi)vdx\\
\geq&(2^*-2)t_n(v)^2\int_\Omega|t_n(v) v+\mu_n\varphi|^{2^*-2}v^2dx-\lambda_0\mu_nt_n(v)\int_\Omega\varphi vdx\\
&-\int_\Omega|t_n(v) v+\mu_n\varphi|^{2^*-2}t_n(v)\mu_n\varphi vdx\\
\geq&Ct_n(v)^{2^*}\|v\|_{2^*}^{2^*}-t_n(v) C_1\lambda_0\mu_n\|v\|_{2^*}-C_1\mu_n t_n(v)^{2^*-1}\|v\|_{2^*}^{2^*-1}-C_1\mu_n^{2^*-2}t_n(v)^2\|v\|_{2^*}^2\\
&-C_1\mu_n^{2^*-1}t_n(v)\|v\|_{2^*}\\
\geq&Ct_n(v)^{2^*}\|v\|_{2^*}^{2^*}-C_1\mu_n^{2^*}-C_1(\lambda_0\mu_n)^{\frac{2^*}{2^*-1}}.
\endaligned
\end{equation}
Since $t_n(v)\geq t_0(v)$, $(t_n(v)\|v\|_{2^*})^{2^*}\geq C((1-\frac{\lambda_0}{\lambda_1})-C_1\mu_n^{2^*-2})^{\frac{2^*}{2^*-2}}$, which implies
$$
0\geq C((1-\frac{\lambda_0}{\lambda_1})-C_1\mu_n^{2^*-2})^{\frac{2^*}{2^*-2}}-C_1\mu_n^{2^*}-C_1(\lambda_0\mu_n)^{\frac{2^*}{2^*-1}}
=C(1-\frac{\lambda_0}{\lambda_1}+o_n(1))^{\frac{2^*}{2^*-2}}+o_n(1),
$$
which is a contradiction for $n$ large enough since $\lambda_0<\lambda_1$.

{\bf Claim~3}\quad $T_{\lambda,\mu}'(t)\to-\infty$ as $t\to+\infty$ for every $\lambda>0$ and $\mu>0$.

Indeed,
$$
\aligned
T_{\lambda,\mu}'(t)&=t(\|v\|-\lambda\|v\|_2^2)-\lambda\mu\int_\Omega\varphi vdx-\int_\Omega|tv+\mu\varphi|^{2^*-2}(tv+\mu\varphi)vdx\\
=&t(\|v\|-\lambda\|v\|_2^2)-\lambda\mu\int_\Omega\varphi vdx-t^{2^*-1}(\int_\Omega|v|^{2^*}dx+o_t(1)).
\endaligned
$$
Since $v\not=0$, $T_{\lambda,\mu}'(t)\to-\infty$ as $t\to+\infty$.

By Claims~1 and 3, for every $\lambda\in(0, \lambda_1)$ and $\mu\in(0, \mu_\lambda^2)$, there exists $t_{\lambda,\mu}^-(v)>t_0(v)$ such that $t_{\lambda,\mu}^-(v)v\in\mathcal{N}_{\lambda,\mu}^-$ for every $v\in\h\backslash\{0\}$.  Thanks to Claim~2, $t_{\lambda,\mu}^-(v)$ is uniqueness when $\lambda\in(0, \lambda_1)$ and $\mu\in(0, \mu_\lambda^3)$.  Moreover, if $\lambda\mu\int_\Omega\varphi vdx+\mu^{2^*-1}\int_\Omega\varphi^{2^*-1}vdx>0$, then by Claim~1, for every $\lambda\in(0, \lambda_1)$ and $\mu\in(0, \mu_\lambda^2)$, there exists unique $t_{\lambda,\mu}^+(v)\in(0, t_0(v))$ such that $t_{\lambda,\mu}^+(v)v\in\mathcal{N}_{\lambda,\mu}^+$ and $T_{\lambda,\mu}(t_{\lambda,\mu}^+(v))=\min_{0\leq t\leq t_{\lambda,\mu}^-(v)}T_{\lambda, \mu}(t)$.
If there exist $\lambda_0\in(0, \lambda_1)$ and $\mu_0\in(0, \mu_{\lambda_0}^3)$ such that $\mathcal{N}_{\lambda_0,\mu_0}^0\not=\emptyset$, then there exists $v_0\in\mathcal{N}_{\lambda_0,\mu_0}^0$.  This means that $T_{\lambda_0,\mu_0,v_0}'(1)=0$ and $T_{\lambda_0,\mu_0,v_0}''(1)=0$.  By Claim~2, we have $1<t_0(v_0)$.  This, together with the choice of $t_0(v_0)$, implies $T_{\lambda_0,\mu_0,v_0}''(1)>0$, which is a contradiction.  Therefore, $\mathcal{N}_{\lambda,\mu}^0=\emptyset$ for every $\lambda\in(0, \lambda_1)$ and $\mu\in(0, \mu_\lambda^3)$.  If there exist $\lambda_0\in(0, \lambda_1)$ and a sequence $\mu_n\to0$ as $n\to\infty$ such that dist$(\mathcal{N}_{\lambda_0,\mu_n}^+, \mathcal{N}_{\lambda_0,\mu_n}^-)\to0$ as $n\to\infty$, then there exist $v_n\in\mathcal{N}_{\lambda_0,\mu_n}^+$ and $z_n\in\mathcal{N}_{\lambda_0,\mu_n}^-$ such that $\|v_n-z_n\|\to0$ as $n\to\infty$.  Since $v_n\in\mathcal{N}_{\lambda_0,\mu_n}^+$, $T_{\lambda_0,\mu_n,v_n}'(1)=0$ and $T_{\lambda_0,\mu_n,v_n}''(1)>0$ for all $n\in\bbn$.  Similar to \eqref{eq999}, we obtain
\begin{equation}\label{eq008}
\|v_n\|< C\big((\frac{\lambda_1}{\lambda_1-\lambda_0})((\lambda_0\mu_n)^{\frac{2^*}{2^*-1}}+\mu_n^{2^*})\big)^{\frac12}\quad\text{for all }n\in\bbn.
\end{equation}
On the other hand, since $z_n\in\mathcal{N}_{\lambda_0,\mu_n}^-$, $T_{\lambda_0,\mu_n,z_n}'(1)=0$ and $T_{\lambda_0,\mu_n,z_n}''(1)<0$ for all $n\in\bbn$.  This implies $1>t_0(z_n)$ for all $n\in\bbn$.  Thus,
\begin{equation}\label{eq007}
\|z_n\|>C\|t_0(z_n)z_n\|_{2^*}\geq(C(1-\frac{\lambda_0}{\lambda_1})-C_1\mu_n^{2^*-2})^{\frac{1}{2^*-2}}\quad\text{for all }n\in\bbn.
\end{equation}
Combining \eqref{eq008} and \eqref{eq007}, we have
$$
o_n(1)=\|v_n-z_n\|\geq\|z_n\|-\|v_n\|=C(1-\frac{\lambda_0}{\lambda_1}+o_n(1))^{\frac{1}{2^*-2}}+o_n(1),
$$
which is a contradiction for $n$ large enough.  Thus, for every $\lambda\in(0, \lambda_1)$, there exists $\mu_\lambda^4\leq\mu_\lambda^3$ such that $\text{dist}(\mathcal{N}_{\lambda,\mu}^-,\mathcal{N}_{\lambda,\mu}^+):=\inf\{\|v-z\|:v\in\mathcal{N}_{\lambda,\mu}^-, z\in\mathcal{N}_{\lambda,\mu}^+\}>C(1-\frac{\lambda}{\lambda_1})^{\frac{1}{2^*-2}}>0$.  We complete the proof of this lemma by taking $\mu_\lambda^0=\mu_\lambda^4$.
\end{prooff}

\begin{Lem}\label{lem002}
For every $\lambda\in(0, \lambda_1)$ and $\mu\in(0, \mu_\lambda^0)$, $\inf_{\mathcal{N}_{\lambda,\mu}}I_{\lambda,\mu}(v)>-\infty$.
\end{Lem}
\begin{prooff}
By Lemma~\ref{lem001}, for every $\lambda\in(0, \lambda_1)$ and $\mu\in(0, \mu_\lambda^0)$, $\mathcal{N}_{\lambda,\mu}\not=\emptyset$.  Let $v\in\mathcal{N}_{\lambda,\mu}$, then, by the Young inequality, we have
\begin{equation}\label{eq006}
\aligned
I_{\lambda, \mu}(v)=&I_{\lambda, \mu}(v)-\frac12I_{\lambda, \mu}'(v)v\\
=&-\frac{\lambda\mu}{2}\int_\Omega\varphi vdx-\frac{\lambda}{2}\|\mu\varphi\|_2^2+\frac12\int_\Omega|v+\mu\varphi|^{2^*-2}(v+\mu\varphi)vdx
-\frac{1}{2^*}\|v+\mu\varphi\|_{2^*}^{2^*}\\
\geq&(\frac12-\frac{1}{2^*})\|v+\mu\varphi\|_{2^*}^{2^*}-\frac12\|v+\mu\varphi\|_{2^*}^{2^*-1}\|\mu\varphi\|_{2^*}
-\frac{\lambda}{2}\|\mu\varphi\|_2^2-\frac{\lambda}{2}\|\mu\varphi\|_2\|v\|_2\\
\geq&C\|v\|_{2^*}^{2^*}-C(\lambda,\mu).
\endaligned
\end{equation}
This means that $\inf_{\mathcal{N}_{\lambda,\mu}}I_{\lambda,\mu}(v)>-C(\lambda,\mu)$.
\end{prooff}

By Lemma~\ref{lem002}, we can define $m_{\lambda,\mu}^\pm=\inf_{\mathcal{N}_{\lambda,\mu}^\pm}I_{\lambda,\mu}(v)$ when $\lambda\in(0, \lambda_1)$ and $\mu\in(0, \mu_\lambda^0)$.  Moreover, by $(4)$ of Lemma~\ref{lem005}, if $v_{\lambda,\mu}^\pm\in\mathcal{N}_{\lambda,\mu}^\pm$ satisfy $I_{\lambda,\mu}(v_{\lambda,\mu}^\pm)=m_{\lambda,\mu}^\pm$, then $v_{\lambda,\mu}^\pm$ are also local minimizers of $I_{\lambda,\mu}$ on $\mathcal{N}_{\lambda,\mu}$.  Next lemma gives an estimate of $m_{\lambda,\mu}^\pm$.
\begin{Lem}\label{lem003}
For every $\lambda\in(0, \lambda_1)$, there exists $\mu_\lambda^{*}\in(0, \mu_\lambda^0]$ such that
$$
m_{\lambda,\mu}^+\leq I_{\lambda,\mu}(0)<0<m_{\lambda,\mu}^-\quad\text{for }\mu\in(0, \mu_\lambda^*).
$$
\end{Lem}
\begin{prooff}
We first estimate $m_{\lambda,\mu}^+$.  Indeed, $T_{\lambda,\mu,v}(t_{\lambda,\mu}^+(v))=\min_{0\leq t\leq t_{\lambda,\mu}^-(v)}T_{\lambda,\mu,v}(t)\leq T_{\lambda,\mu,v}(0)$
for $\lambda\in(0, \lambda_1)$, $\mu\in(0, \mu_\lambda^0)$ and $v\in\h\backslash\{0\}$ by Lemma~\ref{lem001}.  So
$$
m_{\lambda,\mu}^+\leq I_{\lambda,\mu}(t_{\lambda,\mu}^+(v)v)=T_{\lambda,\mu,v}(t_{\lambda,\mu}^+(v))\leq T_{\lambda,\mu,v}(0)=I_{\lambda,\mu}(0)<0.
$$
Next, we estimate $m_{\lambda,\mu}^-$.  Assume $v\in\mathcal{N}_{\lambda,\mu}^-$, similar to \eqref{eq006}, we have that
$$
I_{\lambda, \mu}(v)\geq C\|v\|_{2^*}^{2^*}-C_1\mu^{2^*}-C_1\lambda\mu^2-C_1(\lambda\mu)^{2^*/(2^*-1)}.
$$
Since $v\in\mathcal{N}_{\lambda,\mu}^-$, similar to \eqref{eq007}, we have $\|v\|_{2^*}\geq (C(1-\frac{\lambda}{\lambda_1})-C_1\mu^{2^*-2}))^{1/(2^*-2)}$.  Thus, for every $\lambda\in(0, \lambda_1)$, there exists $\mu_\lambda^*\in(0, \mu_\lambda^0]$ such that $I_{\lambda, \mu}(v)\geq C(1-\frac{\lambda}{\lambda_1})^{2^*/(2^*-2)}>0$ for $v\in\mathcal{N}_{\lambda,\mu}^-$, which implies $m_{\lambda,\mu}^->0$.
\end{prooff}

\begin{Rmk}\label{rmk001}
Combining Lemmas~\ref{lem001} and \ref{lem003}, we can conclude that for every $\lambda\in(0, \lambda_1)$ and $\mu\in(0, \mu_\lambda^*)$, $I_{\lambda,\mu}(t_{\lambda,\mu}^-(v)v)=\max_{t\geq0}I_{\lambda,\mu}(tv)$ for every $v\in\h\backslash\{0\}$.
\end{Rmk}

We close this section by establishing a local compactness lemma.
\begin{Lem}\label{lem004}
Assume $I_{\lambda,\mu}(v_n)=c+o_n(1)$ and $I_{\lambda,\mu}'(v_n)=o_n(1)$.  If $c<\frac1N S^{N/2}+m_{\lambda,\mu}^+$, then there exists $v_0\in\h$ such that $v_n\to v_0$ in $\h$ up to a subsequence.
\end{Lem}
\begin{prooff}
Similar to \eqref{eq006}, we know that $\{v_n\}$ is bounded in $L^{2^*}(\Omega)$.  Since
$$
\aligned
&I_{\lambda,\mu}(v_n)-\frac{1}{2^*}I_{\lambda,\mu}'(v_n)v_n\\
=&(\frac12-\frac1{2^*})(\|v_n\|^2-\lambda\|v_n\|_2^2)
-(1-\frac1{2^*})\lambda\int_\Omega v_n\mu\varphi dx-\frac\lambda2\|\mu\varphi\|_2^2\\
&-\frac1{2^*}\int_\Omega|v_n+\mu\varphi|^{2^*-2}(v_n+\mu\varphi)\mu\varphi dx\\
\geq&(\frac12-\frac1{2^*})(1-\frac{\lambda}{\lambda_1})\|v_n\|^2-C\lambda\mu\|v_n\|_{2^*}-C\lambda\mu^2-C\mu^{2^*}
-C\mu\|v_n\|_{2^*}^{2^*-1},
\endaligned
$$
$\{v_n\}$ is bounded in $\h$.  Without loss of generality, we assume $v_n\rightharpoonup v_0$ in $\h$ as $n\to\infty$.  Let $w_n=v_n-v_0$ then $w_n\rightharpoonup0$ in $\h$ as $n\to\infty$.  Note that, by the Br\'ezis--Lieb Lemma and $v_n\to v_0$ in $L^2(\Omega)$, we have
\begin{equation}\label{eq009}
\aligned
I_{\lambda,\mu}(w_n)=&\frac12\|w_n\|^2-\frac{\lambda}2\|w_n+\mu\varphi\|_2^2-\frac1{2^*}\|w_n+\mu\varphi\|_{2^*}^{2^*}\\
=&I_{\lambda,\mu}(v_n)-I_{\lambda,\mu}(v_0)+I_{\lambda,\mu}(0)+o_n(1)\\
&-\frac{\lambda}{2}\int_\Omega
(|v_n-v_0+\mu\varphi|^2+|v_0+\mu\varphi|^2-|v_n+\mu\varphi|^2-|\mu\varphi|^2)dx\\
&-\frac1{2^*}\int_{\Omega}(|v_n-v_0+\mu\varphi|^{2^*}+|v_0+\mu\varphi|^{2^*}-|v_n+\mu\varphi|^{2^*}-|\mu\varphi|^{2^*})dx\\
=&I_{\lambda,\mu}(v_n)-I_{\lambda,\mu}(v_0)+I_{\lambda,\mu}(0)+o_n(1)\\
&-\frac1{2^*}\int_{\Omega}(|w_n+\mu\varphi|^{2^*}-|w_n|^{2^*}-|\mu\varphi|^{2^*})dx\\
&-\frac1{2^*}\int_{\Omega}(|v_0+\mu\varphi|^{2^*}+|w_n|^{2^*}-|w_n+v_0+\mu\varphi|^{2^*})dx\\
=&I_{\lambda,\mu}(v_n)-I_{\lambda,\mu}(v_0)+I_{\lambda,\mu}(0)+o_n(1).
\endaligned
\end{equation}
On the other hand, also by the Br\'ezis--Lieb Lemma, we have
\begin{equation}\label{eq900}
\aligned
&I_{\lambda,\mu}'(w_n)w_n\\
=&I_{\lambda,\mu}'(v_n)v_n-I_{\lambda,\mu}'(v_0)v_0+o_n(1)-\int_\Omega(|w_n+\mu\varphi|^{2^*-2}(w_n+\mu\varphi)w_ndx\\
&+\int_\Omega(|v_n+\mu\varphi|^{2^*-2}(v_n+\mu\varphi)v_n-|v_0+\mu\varphi|^{2^*-2}(v_0+\mu\varphi)v_0)dx\\
=&I_{\lambda,\mu}'(v_n)v_n-I_{\lambda,\mu}'(v_0)v_0-\int_\Omega(|v_n-v_0+\mu\varphi|^{2^*}+|v_0+\mu\varphi|^{2^*}-|v_n+\mu\varphi|^{2^*}-|\mu\varphi|^{2^*})dx\\
&+o_n(1)+\int_\Omega(|w_n+\mu\varphi|^{2^*-2}(w_n+\mu\varphi)-|\mu\varphi|^{2^*-1}-|v_n+\mu\varphi|^{2^*-2}(v_n+\mu\varphi))\mu\varphi dx\\
&+\int_\Omega|v_0+\mu\varphi|^{2^*-2}(v_0+\mu\varphi)\mu\varphi dx\\
=&I_{\lambda,\mu}'(v_n)v_n-I_{\lambda,\mu}'(v_0)v_0+o_n(1).
\endaligned
\end{equation}
Since $v_n\rightharpoonup v_0$ and $I_{\lambda,\mu}'(v_n)=o_n(1)$, $I_{\lambda,\mu}'(v_0)=0$.  This means that $I_{\lambda,\mu}'(w_n)w_n=o_n(1)$.  Assume $\|w_n\|^2=b+o_n(1)$.  Then, by $I_{\lambda,\mu}'(w_n)w_n=o_n(1)$, we know that $\|w_n\|_{2^*}^{2^*}=b+o_n(1)$.  This, together with the Br\'ezis--Lieb Lemma, implies that
\begin{equation}\label{eq010}
I_{\lambda,\mu}(w_n)=\frac1N b+o_n(1)+I_{\lambda,\mu}(0).
\end{equation}
Note that $v_0\not=0$ and $I_{\lambda,\mu}'(v_0)=0$ for $\mu>0$.  Therefore, $I_{\lambda,\mu}(v_0)\geq m_{\lambda,\mu}^+$ by (3) of Lemma~\ref{lem001} and Lemma~\ref{lem003}.  Thus, combining \eqref{eq009} and \eqref{eq010}, we conclude that $b<S^{N/2}$.  On the other hand, if $b\not=0$ then $b\geq S^{N/2}$, since
$$
b^{2^*/2}S^{-2^*/2}+o_n(1)=\|w_n\|^{2^*}S^{-2^*/2}\geq\|w_n\|_{2^*}^{2^*}=b+o_n(1).
$$
Therefore, we have $v_n\to v_0$ in $\h$.
\end{prooff}

\section{Two solutions of $(\mathcal{Q}_{\lambda,\mu})$}
In this section, we will obtain two solutions of $(\mathcal{Q}_{\lambda,\mu})$ by using the method of fibering maps.  We start by
\begin{Lem}\label{lem006}
For every $\lambda\in(0, \lambda_1)$ and $\mu\in(0, \mu_\lambda^{*})$,  $(\mathcal{Q}_{\lambda,\mu})$ has a positive solution $v_{\lambda,\mu}^+\in\mathcal{N}_{\lambda,\mu}^+$ with $I_{\lambda,\mu}(v_{\lambda,\mu}^+)=m_{\lambda,\mu}^+$.
\end{Lem}
\begin{prooff}
By Lemma~\ref{lem001}, $\mathcal{N}_{\lambda,\mu}^0=\emptyset$ for $\lambda\in(0, \lambda_1)$ and $\mu\in(0, \mu_\lambda^*)$.  Thus, using the Ekeland principle in a standard way (cf. Sun and Li \cite{SL08}), we can obtain a sequence $\{v_{\lambda,\mu}^n\}\subset\mathcal{N}_{\lambda,\mu}^+$ satisfying $I_{\lambda,\mu}(v_{\lambda,\mu}^n)\to m_{\lambda,\mu}^+$ and $I_{\lambda,\mu}'(v_{\lambda,\mu}^n)\to0$ as $n\to\infty$.  Thanks to Lemma~\ref{lem004}, there exists $v_{\lambda,\mu}^+$ such that $v_{\lambda,\mu}^n\to v_{\lambda,\mu}^+$ in $\h$ as $n\to\infty$, which implies $I_{\lambda,\mu}(v_{\lambda,\mu}^+)=m_{\lambda,\mu}^+$.  Without loss of generality, we can choose $v_{\lambda,\mu}^+\geq0$.   Indeed, if $(v_{\lambda,\mu}^+)_-:=\min\{v_{\lambda,\mu}^+, 0\}\not=0$, then $|\mu\varphi+v_{\lambda,\mu}^+|\leq\mu\varphi+|v_{\lambda,\mu}^+|$ since $\varphi\geq0$.  Therefore,
\begin{equation}\label{eq014}
\aligned
&I_{\lambda,\mu}(v_{\lambda,\mu}^+)\\
=&\frac12\|v_{\lambda,\mu}^+\|^2-\frac\lambda2\|v_{\lambda,\mu}^++\mu\varphi\|_2^2-
\frac1{2^*}\|v_{\lambda,\mu}^++\mu\varphi\|_{2^*}^{2^*}\\
\geq&\frac12\big\||v_{\lambda,\mu}^+|\big\|^2-\frac\lambda2\big\||v_{\lambda,\mu}^+|+\mu\varphi\big\|_2^2-
\frac1{2^*}\big\||v_{\lambda,\mu}^+|+\mu\varphi\big\|_{2^*}^{2^*}\\
=&I_{\lambda,\mu}(|v_{\lambda,\mu}^+|).
\endaligned
\end{equation}
On the other hand, by Lemma~\ref{lem001}, there exist $0<t_{\lambda,\mu}^+(|v_{\lambda,\mu}^+|)<t_{\lambda,\mu}^-(|v_{\lambda,\mu}^+|)<+\infty$ such that $t_{\lambda,\mu}^+(|v_{\lambda,\mu}^+|)|v_{\lambda,\mu}^+|\in\mathcal{N}_{\lambda,\mu}^+$ and $t_{\lambda,\mu}^-(|v_{\lambda,\mu}^+|)|v_{\lambda,\mu}^+|\in\mathcal{N}_{\lambda,\mu}^-$.  There are two cases may occur:
\begin{itemize}
\item[$(1)$] $t_{\lambda,\mu}^-(|v_{\lambda,\mu}^+|)<1$;
\item[$(2)$] $1\leq t_{\lambda,\mu}^-(|v_{\lambda,\mu}^+|)$.
\end{itemize}
If case~(1) holds, then $1>t_{\lambda,\mu}^-(|v_{\lambda,\mu}^+|)>t_0(|v_{\lambda,\mu}^+|)$ since $t_{\lambda,\mu}^-(|v_{\lambda,\mu}^+|)|v_{\lambda,\mu}^+|\in\mathcal{N}_{\lambda,\mu}^-$.  It follows that
$$
\big\||v_{\lambda,\mu}^+|\big\|_{2^*}>\big\|t_{\lambda,\mu}^-(|v_{\lambda,\mu}^+|)|v_{\lambda,\mu}^+|\big\|_{2^*}>
\big\|t_0(|v_{\lambda,\mu}^+|)|v_{\lambda,\mu}^+|\big\|_{2^*}.
$$
On the other hand, $1<t_0(v_{\lambda,\mu}^+)$ since $v_{\lambda,\mu}^+\in\mathcal{N}_{\lambda,\mu}^+$.  It follows that
$$
\big\||v_{\lambda,\mu}^+|\big\|_{2^*}=\|v_{\lambda,\mu}^+\|_{2^*}<\|t_0(v_{\lambda,\mu}^+)v_{\lambda,\mu}^+\|_{2^*}.
$$
Thus, $t_0(v_{\lambda,\mu}^+)>t_0(|v_{\lambda,\mu}^+|)$.  Recall the choice of $t_0(v)$, we know that $t_0(v_{\lambda,\mu}^+)=t_0(|v_{\lambda,\mu}^+|)$, which is a contradiction.  Therefore, case $(2)$ must happen.
In this case, by $(2)$ of Lemma~\ref{lem001}, we have
$$
I_{\lambda,\mu}(t_{\lambda,\mu}^+(|v_{\lambda,\mu}^+|)|v_{\lambda,\mu}^+|)=\min_{0\leq t\leq t_{\lambda,\mu}^-(|v_{\lambda,\mu}^+|)}I_{\lambda,\mu}(t|v_{\lambda,\mu}^+|)\leq I_{\lambda,\mu}(|v_{\lambda,\mu}^+|).
$$
This, together with \eqref{eq014}, gives that
$$
m_{\lambda,\mu}^+=I_{\lambda,\mu}(v_{\lambda,\mu}^+)\geq I_{\lambda,\mu}(|v_{\lambda,\mu}^+|)\geq I_{\lambda,\mu}(t_{\lambda,\mu}^+(|v_{\lambda,\mu}^+|)|v_{\lambda,\mu}^+|)\geq m_{\lambda,\mu}^+,
$$
that is, $I_{\lambda,\mu}(t_{\lambda,\mu}^+(|v_{\lambda,\mu}^+|)|v_{\lambda,\mu}^+|)=m_{\lambda,\mu}^+$.  Thanks to $(3)$ and $(4)$ of Lemma~\ref{lem001}, $t_{\lambda,\mu}^+(|v_{\lambda,\mu}^+|)|v_{\lambda,\mu}^+|$ is also a local minimum of $I_{\lambda,\mu}$ on $\mathcal{N}_{\lambda,\mu}$ for $\lambda\in(0, \lambda_1)$ and $\mu\in(0, \mu_\lambda^*)$.  Similar to \cite[Theorem~2.3]{BZ03}, we can show that $I_{\lambda,\mu}'(t_{\lambda,\mu}^+(|v_{\lambda,\mu}^+|)|v_{\lambda,\mu}^+|)=0$.  Therefore, by the maximum principle, there exists $v_{\lambda,\mu}^+>0$ such that $I_{\lambda,\mu}(v_{\lambda,\mu}^+)=m_{\lambda,\mu}^+$ and $I_{\lambda,\mu}'(v_{\lambda,\mu}^+)=0$ for $\lambda\in(0, \lambda_1)$ and $\mu\in(0, \mu_\lambda^*)$.
\end{prooff}

\begin{Lem}\label{lem020}
For every $\lambda\in(0, \lambda_1)$ and $\mu\in(0, \mu_\lambda^*)$, $t_{\lambda,\mu}^-(u)$ is continuous for $u\in\h\backslash\{0\}$.  Moreover, $\mathcal{N}_{\lambda,\mu}^-=\{u\in\h\backslash\{0\}: t_{\lambda,\mu}^-(\frac{u}{\|u\|})\frac{1}{\|u\|}=1\}$.
\end{Lem}
\begin{prooff}
Since $t_{\lambda,\mu}^-(u)$ is unique for every $u\in\h\backslash\{0\}$, $t_{\lambda,\mu}^-(u)$ is continuous for $u\in\h\backslash\{0\}$.  Let $v=\frac{u}{\|u\|}$, then by Lemma~\ref{lem001}, there exists $t_{\lambda,\mu}^-(v)>0$ such that $t_{\lambda,\mu}^-(v)v\in\mathcal{N}_{\lambda,\mu}^-$, that is, $t_{\lambda,\mu}^-(\frac{u}{\|u\|})\frac{u}{\|u\|}\in\mathcal{N}_{\lambda,\mu}^-$.  If $u\in\mathcal{N}_{\lambda,\mu}^-$, then by the uniqueness of $t_{\lambda,\mu}^-(u)$, we must have $t_{\lambda,\mu}^-(\frac{u}{\|u\|})\frac{1}{\|u\|}=1$.  Therefore, $\mathcal{N}_{\lambda,\mu}^-\subset\{u\in\h\backslash\{0\}: t_{\lambda,\mu}^-(\frac{u}{\|u\|})\frac{1}{\|u\|}=1\}$.  On the other hand, if $t_{\lambda,\mu}^-(\frac{u}{\|u\|})\frac{1}{\|u\|}=1$, then also by Lemma~\ref{lem001}, $u=t_{\lambda,\mu}^-(\frac{u}{\|u\|})\frac{u}{\|u\|}\in\mathcal{N}_{\lambda,\mu}^-$, which implies that $\mathcal{N}_{\lambda,\mu}^-=\{u\in\h\backslash\{0\}: t_{\lambda,\mu}^-(\frac{u}{\|u\|})\frac{1}{\|u\|}=1\}$.
\end{prooff}

\begin{Lem}\label{lem021}
For every $\lambda\in(0, \lambda_1)$ and $\mu\in(0, \mu_\lambda^*)$, $m_{\lambda,\mu}^-<m_{\lambda,\mu}^++\frac1N S^{N/2}$.
\end{Lem}
\begin{prooff}
Let $w_{\lambda,\mu,t}=v_{\lambda,\mu}^++tU_\lambda$, where $U_\lambda$ is the ground state of $(\mathcal{P}_{\lambda, 0})$.  For the sake of clarity, we divide the proof into three steps.

{\bf Step.~1}\quad For every $\lambda\in(0, \lambda_1)$ and $\mu\in(0, \mu_\lambda^*)$, there exists $C_{\lambda,\mu}>0$
such that $t_{\lambda,\mu}^-(\widetilde{w}_{\lambda,\mu,t})<C_{\lambda,\mu}$ for all $t\geq0$, where $\widetilde{w}_{\lambda,\mu,t}=\frac{w_{\lambda,\mu,t}}{\|w_{\lambda,\mu,t}\|}$.

Indeed, by Lemma~\ref{lem020}, if not, then there exists $\{t_n\}\subset\bbr$ with $t_n\to+\infty$ such that $t_{\lambda,\mu}^-(\widetilde{w}_{\lambda,\mu,t_n})\to+\infty$.  By the Lebesgue dominated convergence theorem,
$$
\int_\Omega|\widetilde{w}_{\lambda,\mu,t_n}|^{2^*}dx=\int_\Omega\frac{|v_{\lambda,\mu}^++t_nU_\lambda|^{2^*}}
{\|v_{\lambda,\mu}^++t_nU_\lambda\|^{2^*}}dx=\frac{\|U_\lambda\|^{2^*}_{2^*}}{\|U_\lambda\|^{2^*}}+o_n(1).
$$
This, together with Lemma~\ref{lem003} and the Lebesgue dominated convergence theorem, implies
$$
\aligned
0<&I_{\lambda,\mu}(t_{\lambda,\mu}^-(\widetilde{w}_{\lambda,\mu,t_n})\widetilde{w}_{\lambda,\mu,t_n})\\
=&\frac{(t_{\lambda,\mu}^-(\widetilde{w}_{\lambda,\mu,t_n}))^2}{2}\|\widetilde{w}_{\lambda,\mu,t_n}\|^2
-\frac\lambda2\int_\Omega (t_{\lambda,\mu}^-(\widetilde{w}_{\lambda,\mu,t_n})\widetilde{w}_{\lambda,\mu,t_n}+\mu\varphi)^2dx\\
&-\frac{1}{2^*}\int_\Omega(t_{\lambda,\mu}^-(\widetilde{w}_{\lambda,\mu,t_n})\widetilde{w}_{\lambda,\mu,t_n}+\mu\varphi)^{2^*}dx\\
=&\frac{(t_{\lambda,\mu}^-(\widetilde{w}_{\lambda,\mu,t_n}))^2}{2}(\|\widetilde{w}_{\lambda,\mu,t_n}\|^2-\lambda\|\widetilde{w}_{\lambda,\mu,t_n}\|_2^2+o_n(1))
-\frac{(t_{\lambda,\mu}^-(\widetilde{w}_{\lambda,\mu,t_n}))^{2^*}}{2^*}(\|\widetilde{w}_{\lambda,\mu,t_n}\|_{2^*}^{2^*}+o_n(1))\\
\leq&\frac{(t_{\lambda,\mu}^-(\widetilde{w}_{\lambda,\mu,t_n}))^2}{2}(1+o_n(1))-\frac{(t_{\lambda,\mu}^-(\widetilde{w}_{\lambda,\mu,t_n})^{2^*}}{2^*}
(\frac{\|U_\lambda\|^{2^*}_{2^*}}{\|U_\lambda\|^{2^*}}+o_n(1))\to-\infty\quad\text{as }n\to\infty,
\endaligned
$$
a contradiction.

{\bf Step.~2}\quad For every $\lambda\in(0, \lambda_1)$ and $\mu\in(0, \mu_\lambda^*)$, there exists
$t_{\lambda,\mu}>0$ such that $v_{\lambda,\mu}^++t_{\lambda,\mu}U_\lambda\in\mathcal{N}_{\lambda,\mu}^-$.

Let $\mathcal{A}_{\lambda,\mu}^-:=\{u\in\h\backslash\{0\}: t_{\lambda,\mu}^-(\frac{u}{\|u\|})\frac{1}{\|u\|}<1\}$ and $\mathcal{A}_{\lambda,\mu}^+:=\{u\in\h\backslash\{0\}: t_{\lambda,\mu}^-(\frac{u}{\|u\|})\frac{1}{\|u\|}>1\}\cup\{0\}$.  By Lemma~\ref{lem020}, we know that $\h=\mathcal{A}_{\lambda,\mu}^-\cup\mathcal{A}_{\lambda,\mu}^+\cup\mathcal{N}_{\lambda,\mu}^-$.  By Lemma~\ref{lem001}, it is easy to see that $\mathcal{N}_{\lambda,\mu}^+\subset \mathcal{A}_{\lambda,\mu}^+$. In particular, $v_{\lambda,\mu}^+\in \mathcal{A}_{\lambda,\mu}^+$. On the other hand, let $t_{\lambda,\mu}^1=\frac{C_{\lambda,\mu}+1+\|v_{\lambda,\mu}^+\|}{\|U_\lambda\|}$, where $C_{\lambda,\mu}$ is given in Step.~1.  Then by Step.~1,
$$
\|v_{\lambda,\mu}^++t_{\lambda,\mu}^1U_\lambda\|^2\geq(t_{\lambda,\mu}^1\|U_\lambda\|-\|v_{\lambda, \mu}^+\|)^2= (C_{\lambda,\mu}+1)^2>\bigg(t(\frac{v_{\lambda,\mu}^++t_{\lambda,\mu}^1U_\lambda}{\|v_{\lambda,\mu}^++t_{\lambda,\mu}^1U_\lambda\|})^-\bigg)^2,
$$
which implies that $t_{\lambda,\mu}^-(\frac{v_{\lambda,\mu}^++t_{\lambda,\mu}^1U_\lambda}{\|v_{\lambda,\mu}^++t_{\lambda,\mu}^1U_\lambda\|})\frac{1}{\|v_{\lambda,\mu}^++t_{\lambda,\mu}^1U_\lambda\|}<1$, that is, $v_{\lambda,\mu}^++t_{\lambda,\mu}^1U_\lambda\in\mathcal{A}_{\lambda,\mu}^-$.  By the continuity of $t_{\lambda,\mu}^-(u)$ for $u$, there exists $t_{\lambda,\mu}\in(0, t_{\lambda,\mu}^1)$ such that $v_{\lambda,\mu}^++t_{\lambda,\mu}U_\lambda\in\{u\in\h\backslash\{0\}: t_{\lambda,\mu}^-(\frac{u}{\|u\|})\frac{1}{\|u\|}=1\}$. Thanks to Lemma~\ref{lem020}, $v_{\lambda,\mu}^++t_{\lambda,\mu}U_\lambda\in\mathcal{N}_{\lambda,\mu}^-$.

{\bf Step.~3}\quad $I_{\lambda,\mu}(v_{\lambda,\mu}^++t_{\lambda,\mu}U_\lambda)<\frac1N S^{N/2}+m_{\lambda,\mu}^+$, where $t_{\lambda,\mu}$ is given in Step.~2.

Since $v_{\lambda,\mu}^+$ is a solution of $(\mathcal{Q}_{\lambda,\mu})$, by a direct calculation, we obtain
$$
\aligned
&I_{\lambda,\mu}(v_{\lambda,\mu}^++t_{\lambda,\mu}U_\lambda)\\
=&I_{\lambda,\mu}(v_{\lambda,\mu}^+)
+I_{\lambda,0}(t_{\lambda,\mu}U_\lambda)\\
&-\frac1{2^*}\int_\Omega((v_{\lambda,\mu}^++t_{\lambda,\mu}U_\lambda+\mu\varphi)^{2^*}
-2^*(v_{\lambda,\mu}^++\mu\varphi)^{2^*-1}t_{\lambda,\mu}U_\lambda-(t_{\lambda,\mu}U_\lambda)^{2^*}-(v_{\lambda,\mu}^++\mu\varphi)^{2^*})dx.
\endaligned
$$
It is well-known that $I_{\lambda,0}(t_{\lambda,\mu}U_\lambda)\leq I_{\lambda, 0}(U_\lambda)<\frac1N S^{N/2}$ (cf. Struwe \cite{S00}).  On the other hand, $(a+b)^p-a^p-b^p-pa^{p-1}b\geq0$ for all $a>0$, $b>0$ and $p>1$.  Therefore, $I_{\lambda,\mu}(v_{\lambda,\mu}^++t_{\lambda,\mu}U_\lambda)<m_{\lambda,\mu}^++\frac1N S^{N/2}$.

Combining Step.~2 and Step.~3, we know that $m_{\lambda,\mu}^-<m_{\lambda,\mu}^++\frac1N S^{N/2}$ for every $\lambda\in(0, \lambda_1)$ and $\mu\in(0, \mu_\lambda^*)$.
\end{prooff}

\begin{Lem}\label{lem007}
For every $\lambda\in(0, \lambda_1)$ and $\mu\in(0, \mu_\lambda^*)$, $(\mathcal{Q}_{\lambda,\mu})$ has a positive solution $v_{\lambda,\mu}^-\in\mathcal{N}_{\lambda,\mu}^-$ with $I_{\lambda,\mu}(v_{\lambda,\mu}^-)=m_{\lambda,\mu}^-$.
\end{Lem}
\begin{prooff}
By Lemma~\ref{lem001}, $\mathcal{N}_{\lambda,\mu}^0=\emptyset$ for $\lambda\in(0, \lambda_1)$ and $\mu\in(0, \mu_\lambda^*)$.  Thus, using the Ekeland principle in a standard way (cf. Sun and Li \cite{SL08}), we can obtain a sequence $\{v_{\lambda,\mu}^n\}\subset\mathcal{N}_{\lambda,\mu}^-$ satisfying $I_{\lambda,\mu}(v_{\lambda,\mu}^n)\to m_{\lambda,\mu}^-$ and $I_{\lambda,\mu}'(v_{\lambda,\mu}^n)\to0$ as $n\to\infty$.  By Lemma~\ref{lem021}, $m_{\lambda,\mu}^-<\frac1N S^{N/2}+m_{\lambda,\mu}^+$.  Thus, by Lemma~\ref{lem004}, there exists $v_{\lambda,\mu}^-\in\h$ such that $v_{\lambda,\mu}^n\to v_{\lambda,\mu}^-$ in $\h$ as $n\to\infty$.  This implies $I_{\lambda,\mu}(v_{\lambda,\mu}^-)=m_{\lambda,\mu}^-$ and $I_{\lambda,\mu}'(v_{\lambda,\mu}^-)=0$.  Without loss of generality, we may choose $v_{\lambda,\mu}^-\geq0$.  Indeed,  by Lemma~\ref{lem001}, there exists $t_{\lambda,\mu}^-(|v_{\lambda,\mu}^-|)>0$ such that $t_{\lambda,\mu}^-(|v_{\lambda,\mu}^-|)|v_{\lambda,\mu}^-|\in\mathcal{N}_{\lambda,\mu}^-$.  Thus, by Remark~\ref{rmk001} and a similar argument of \eqref{eq014}, we have that
$$
m_{\lambda,\mu}^-=I_{\lambda,\mu}(v_{\lambda,\mu}^-)\geq I_{\lambda,\mu}(t_{\lambda,\mu}^-(|v_{\lambda,\mu}^-|)v_{\lambda,\mu}^-)\geq I_{\lambda,\mu}(t_{\lambda,\mu}^-(|v_{\lambda,\mu}^-|)|v_{\lambda,\mu}^-|)\geq m_{\lambda,\mu}^-,
$$
which implies $I_{\lambda,\mu}(t_{\lambda,\mu}^-(|v_{\lambda,\mu}^-|)|v_{\lambda,\mu}^-|)=m_{\lambda,\mu}^-$.  Similar to \cite[Theorem~2.3]{BZ03}, we can show that $I_{\lambda,\mu}'(t_{\lambda,\mu}^-(|v_{\lambda,\mu}^-|)|v_{\lambda,\mu}^-|)=0$ since, by $(4)$ of Lemma~\ref{lem001}, $t_{\lambda,\mu}^-(|v_{\lambda,\mu}^-|)|v_{\lambda,\mu}^-|$ is a local minimum of $I_{\lambda,\mu}$ on $\mathcal{N}_{\lambda,\mu}$.  Therefore, by the maximum principle, there exists $v_{\lambda,\mu}^->0$ such that $I_{\lambda,\mu}(v_{\lambda,\mu}^-)=m_{\lambda,\mu}^-$ and $I_{\lambda,\mu}'(v_{\lambda,\mu}^-)=0$ for $\lambda\in(0, \lambda_1)$ and $\mu\in(0, \mu_\lambda^*)$.
\end{prooff}

Now, we can prove Theorem~\ref{thm003}
\medskip\par\noindent{\bf Proof of Theorem~\ref{thm003}}\quad It follows immediately from Proposition~\ref{prop01}, Lemmas~\ref{lem006} and \ref{lem007}.
\qquad\raisebox{-0.5mm}{\rule{1.5mm}{4mm}}\vspace{6pt}

\section{The third solution $(\mathcal{Q}_{\lambda,\mu})$}
In this section, we will find the third solution of $(\mathcal{Q}_{\lambda,\mu})$ by using the theory of category as in Wu \cite{W10}.  We first recall the definition of Lusternik--Schnirelman category and some useful lemmas.
\begin{Def}
\begin{itemize}
\item[$(i)$] For a topological space $X$, we say a non-empty, closed subset $Y\subset X$ is contractible to a point in $X$ if and only if there exists a continuous mapping
    $$
    \xi:[0, 1]\times Y\to X
    $$
    such that for some $x_0\in X$
    $$
    \xi(0,x)=x\quad\text{for all }x\in Y,
    $$
    and
    $$
    \xi(1,x)=x_0\quad\text{for all }x\in Y.
    $$
\item[$(ii)$] We define
$$
\aligned
\text{cat}(X):=\min&\{k\in\bbn:\text{there exist closed subsets }Y_1,\cdots, Y_k\subset X\\
&\text{ such that }Y_j\text{ is contractible to a point in }X\text{ for all }\\
&j\text{ and }\cup_{j=1}^kY_j=X\}.
\endaligned
$$
\end{itemize}
When there do not exist finitely many closed subsets $Y_1,\cdots, Y_k\subset X$ such that $Y_j$ is contractible to a point in $X$ for all $j$ and $\cup_{j=1}^kY_j=X$, we say cat$(X)=\infty$.
\end{Def}

We also need the following three lemmas.
\begin{Lem}\cite[Theorem~2.3]{A92}\label{lem009}
Suppose that $X$ is a Hilbert manifold and $F\in C^1(X, \bbr)$.  Assume that there are $c_0\in\bbr$ and $k\in\bbn$,
\begin{itemize}
\item[$(i)$] $F(x)$ satisfies the Palais--Smale condition for energy level $c\leq c_0$,
\item[$(ii)$] cat$(\{x\in X:F(x)\leq c_0\})\geq k$.
\end{itemize}
Then $F(x)$ has at least $k$ critical points in $\{x\in X: F(x)\leq c_0\}$.
\end{Lem}
\begin{Lem}\cite[Lemma~2.5]{AT00}\label{lem010}
Let $X$ be a topological space.  Suppose that there are two continuous maps
$$
\Phi: \mathbb{S}^{N-1}\to X,\quad\Psi:X\to \mathbb{S}^{N-1}
$$
such that $\Psi\circ\Phi$ is homotopic to the identity map of $\mathbb{S}^{N-1}$, that is there exists a continuous map $\zeta:[0, 1]\times \mathbb{S}^{N-1}\to \mathbb{S}^{N-1}$ such that
$$
\aligned
&\zeta(0,x)=\Psi\circ\Phi(x)\quad\text{ for each }x\in \mathbb{S}^{N-1},\\
&\zeta(1,x)=x\quad\text{ for each }x\in \mathbb{S}^{N-1}.
\endaligned
$$
Then cat$(X)\geq2$.
\end{Lem}
\begin{Lem}\cite[Lemma~4.4]{W10}\label{lem011}
Assume $\Omega$ satisfies condition $(D)$, then there exists $d_0>0$ such that for $v\in\mathcal{N}_{0,0}$ with $I_{0,0}(v)\leq \frac1N S^{N/2}+d_0$, we have
$$
\int_{\bbrn}\frac{x}{|x|}|\nabla v|^2dx\not=0.
$$
\end{Lem}

The following lemmas are crucial in finding the third solution.
\begin{Lem}\label{lem012}
There exist $\lambda^*\in(0, \lambda_1)$ and $\mu^{**}\in(0, \mu_\lambda^*)$ such that for $\lambda\in(0, \lambda^*)$, $\mu\in(0, \mu^{**})$ and $v\in\mathcal{N}_{\lambda,\mu}^-$ with $I_{\lambda,\mu}(v)<m_{\lambda,\mu}^++\frac1N S^{N/2}$, we have
$$
\int_{\bbrn}\frac{x}{|x|}|\nabla v|^2dx\not=0.
$$
\end{Lem}
\begin{prooff}
Assume $v\in\mathcal{N}_{\lambda,\mu}^-$ and $I_{\lambda,\mu}(v)<m_{\lambda,\mu}^++\frac1N S^{N/2}$.  Similar to \eqref{eq006}, for $\lambda\in(0, \lambda_1)$ and $\mu\in(0, \mu_\lambda^*)$, we have
\begin{equation}\label{eq015}
\frac1N S^{N/2}>I_{\lambda,\mu}(v)-\frac1{2^*}I_{\lambda,\mu}'(v)v\geq C\|v\|_{2^*}^{2^*}-C_1\mu^{2^*}-C_1\lambda\mu^2-C_1(\lambda\mu)^{2^*/(2^*-1)},
\end{equation}
since, by Lemma~\ref{lem003}, $m_{\lambda,\mu}^+<0$.
It is well-known that there exists $t_v>0$ such that $t_vv\in\mathcal{N}_{0,0}$.  Moreover, by Remark~\ref{rmk001}, $I_{\lambda,\mu}(v)\geq I_{\lambda,\mu}(t_vv)$ for $\lambda\in(0, \lambda_1)$ and $\mu\in(0, \mu_\lambda^*)$.  Since $v\in\mathcal{N}_{\lambda,\mu}^-$,  similar to \eqref{eq007}, we have $\|v\|_{2^*}\geq (C(1-\frac{\lambda}{\lambda_1})-C_1\mu^{2^*-2})^{1/(2^*-2)}>0$ for $\lambda\in(0, \lambda_1)$ and $\mu\in(0, \mu^*_\lambda)$.  Thus, there exist $\lambda^0\in(0, \lambda_1)$ and $\mu^*>0$ such that $\|v\|_{2^*}\geq C>0$ for $\lambda\in(0, \lambda^0)$ and $\mu\in(0, \mu^*)$.  This implies that there exists $T>0$ such that $t_v\leq T$ for all $\lambda\in(0, \lambda^0)$, $\mu\in(0, \mu^*)$ and $v\in\mathcal{N}_{\lambda,\mu}^-$.  On the other hand, since $t_vv\in\mathcal{N}_{0,0}$, $\|t_vv\|_{2^*}\geq C>0$.  This, together with \eqref{eq015}, implies that there exists $t_0>0$ such that $t_0\leq t_v$ for all $\lambda\in(0, \lambda^0)$, $\mu\in(0, \mu^*)$ and $v\in\mathcal{N}_{\lambda,\mu}^-$.  Therefore, by the mean value theorem,
$$
\aligned
I_{0,0}(t_vv)=&I_{\lambda,\mu}(t_vv)+\frac\lambda2\|t_vv+\mu\varphi\|_2^2+\frac1{2^*}\int_\Omega(|t_vv+\mu\varphi|^{2^*}-|t_vv|^{2^*})dx\\
\leq&I_{\lambda,\mu}(v)+C(\lambda,\mu)\\
\leq&\frac1N S^{N/2}+C(\lambda,\mu),
\endaligned
$$
where $C(\lambda,\mu)\to0$ as $\lambda\to0$ and $\mu\to0$.
Thus, there exist $\lambda^*\in(0, \lambda^0]$ and $\mu^{**}\in(0, \mu^*]$ such that for $\lambda\in(0, \lambda^*)$, $\mu\in(0, \mu^{**})$ and $v\in\mathcal{N}_{\lambda,\mu}^-$ with $I_{\lambda,\mu}(v)<m_{\lambda,\mu}^++\frac1N S^{N/2}$, we have
$$
I_{0,0}(t_vv)\leq \frac1N S^{N/2}+d_0.
$$
By Lemma~\ref{lem011},
$$
t_v^2\int_{\bbrn}\frac{x}{|x|}|\nabla v|^2dx\not=0.
$$
Since $t_v\geq t_0$ for all $\lambda\in(0, \lambda^*)$, $\mu\in(0, \mu^{**})$ and $v\in\mathcal{N}_{\lambda,\mu}^-$,
$$
\int_{\bbrn}\frac{x}{|x|}|\nabla v|^2dx\not=0
$$
for all $\lambda\in(0, \lambda^*)$, $\mu\in(0, \mu^{**})$ and $v\in\mathcal{N}_{\lambda,\mu}^-$ with $I_{\lambda,\mu}(v)<m_{\lambda,\mu}^++\frac1N S^{N/2}$.
\end{prooff}

Let
$$
U_{\ve,\overrightarrow{y}}(x)=\frac{\phi(x)(N(N-2)\ve)^{(N-2)/2}}{(\ve^2+|x-(1-\ve)\overrightarrow{y}|)^{(N-2)/2}},
$$
where $\phi(x)\in C^\infty_0(\Omega)$ is radially symmetric function such that $0\leq\phi(x)\leq1$ and $\phi(x)\equiv1$ for $2\delta_0\leq|x|\leq\frac1{2\delta_0}$, $\delta_0$ is given in $(D)$, $\overrightarrow{y}\in \mathcal{S}^{N-1}:=\{x\in\bbrn: |x|=1\}$ and $\ve\in(0, 1)$ small enough.  Then we have

\begin{Lem}\label{lem013}
There exist $\ve_0>0$, $\lambda_{**}^1\leq\lambda^*$ and $\mu_{**}^1\leq\mu^{**}$ such that for each $\ve\in(0, \ve_0)$, $\lambda\in(0, \lambda_{**}^1)$ and $\mu\in(0, \mu_{**}^1)$, $v_{\lambda,\mu}^++t_{\lambda,\mu,\ve}U_{\ve,\overrightarrow{y}}(x)\in\mathcal{N}_{\lambda,\mu}^-$ for all $\overrightarrow{y}\in\mathcal{S}^{N-1}$ and some $t_{\lambda,\mu,\ve}>0$.  Moreover, there exist $t_{\lambda,\mu}^*>t_{\lambda,\mu}'>0$, independent of $\ve$, such that $t_{\lambda,\mu,\ve}\in(t_{\lambda,\mu}', t_{\lambda,\mu}^*)$.
\end{Lem}
\begin{prooff}
The proof is similar to the proofs of Step.~1 and Step.~2 in Lemma~\ref{lem021}, so we only sketch it.  By \cite[Lemma~4.2]{HY09}, there exist $t_{\lambda,\mu}^{**}>0$ and $\ve_1>0$ such that for $t\geq t_{\lambda,\mu}^{**}$ and $\ve\in(0, \ve_1)$, we have
$$
\frac{\|v_{\lambda,\mu}^++tU_{\ve,\overrightarrow{y}}\|_{2^*}^{2^*}}{\|v_{\lambda,\mu}^++tU_{\ve,\overrightarrow{y}}\|^{2^*}}
\geq\frac{2^*\|U_{\ve,\overrightarrow{y}}\|_{2^*}^{2^*}}{\|U_{\ve,\overrightarrow{y}}\|^{2^*}-\frac{2^*}{t^{2^*}}\|v_{\lambda,\mu}^+\|^{2^*}}
=\frac{S^{N/2}+o(\ve)}{(S^{N/2})^{2^*/2}+o(\ve)-\frac{2^*}{t^{2^*}}\|v_{\lambda,\mu}^+\|^{2^*}}\geq C>0.
$$
Hence, as in Step.~1 of Lemma~\ref{lem021}, we can obtain $t_{\lambda,\mu}^-(\widetilde{w}_{t,\ve})\leq T$ for some $T>0$ if $t\geq t_{\lambda,\mu}^{**}$ and $\ve\in(0, \ve_1)$, where
$$
\widetilde{w}_{t,\ve}=\frac{v_{\lambda,\mu}^++tU_{\ve,\overrightarrow{y}}}{\|v_{\lambda,\mu}^++tU_{\ve,\overrightarrow{y}}\|}.
$$
Let $t_{\lambda,\mu,\ve}^*=\max\{t_{\lambda,\mu}^{**}, \frac{T+1+\|v_{\lambda,\mu}^+\|}{\|U_{\ve,\overrightarrow{y}}\|}\}$, then also by \cite[Lemma~4.2]{HY09}, there exists $t_{\lambda,\mu}^*>0$ such that $t_{\lambda,\mu,\ve}^*\leq t_{\lambda,\mu}^*$ for $\ve\in(0, \ve_1)$.  Similar to Step.~2 of Lemma~\ref{lem021}, we have
$$
t_{\lambda,\mu}^-(\widetilde{w}_{t_{\lambda,\mu}^*,\ve})\frac{1}{\|v_{\lambda,\mu}^++t_{\lambda,\mu}^*U_{\ve,\overrightarrow{y}}\|}<1.
$$
Thus, as in Step.~2 of Lemma~\ref{lem021}, there exists $t_{\lambda,\mu,\ve}\leq t_{\lambda,\mu}^*$ such that $v_{\lambda,\mu}^++t_{\lambda,\mu,\ve}U_{\ve,\overrightarrow{y}}\in\mathcal{N}_{\lambda,\mu}^-$ for every $\overrightarrow{y}\in\mathcal{S}^{N-1}$.  Moreover, since $v_{\lambda,\mu}^++t_{\lambda,\mu,\ve}U_{\ve,\overrightarrow{y}}\in\mathcal{N}_{\lambda,\mu}^-$, $v_{\lambda,\mu}^+\in\mathcal{N}_{\lambda,\mu}^+$ and $\|U_{\ve,\overrightarrow{y}}\|_{2^*}^{2^*}=S^{N/2}+O(\ve^N)$ (cf. He and Yang \cite{HY09}, Sun and Li \cite{SL08}, Wu \cite{W10}), by a similar argument of \eqref{eq008} and \eqref{eq007}, we have
$$
t_{\lambda,\mu,\ve}\geq\frac{(C(1-\frac{\lambda}{\lambda_1})-C_1\mu^{2^*-2})^{\frac{2^*}{2^*-2}}-C_1\mu^{2^*}-C_1(\lambda\mu)^{\frac{2^*}{2^*-1}}}{S^{N/2}+O(\ve^N)}
$$
Therefore, there exist $\ve_0\leq\ve_1$, $\lambda_{**}^1\leq\lambda^*$ and $\mu_{**}^1\leq\mu^{**}$ such that $t_{\lambda,\mu,\ve}>t_{\lambda,\mu}'>0$ for all $\ve\in(0, \ve_0)$, $\lambda\in(0, \lambda_{**}^1)$ and $\mu\in(0, \mu_{**}^1)$.
\end{prooff}

\begin{Lem}\label{lem014}
Assume $\lambda\in(0, \lambda_{**}^1)$, $\mu\in(0, \mu_{**}^1)$ and $N\geq4$.  Then there exists $\ve^*\leq\ve_0$ such that for $\ve\in(0, \ve^*)$, $v_{\lambda,\mu}^++t_{\lambda,\mu,\ve}U_{\ve,\overrightarrow{y}}\in\mathcal{G}_{\lambda,\mu,\sigma_\ve}$ for all $\overrightarrow{y}\in\mathcal{S}^{N-1}$, where $\mathcal{G}_{\lambda,\mu,\sigma_\ve}:=\{u\in\mathcal{N}_{\lambda,\mu}^-:u\geq0, I_{\lambda,\mu}(u)<\frac1N S^{N/2}+m_{\lambda,\mu}^+-\sigma_\ve\}$ and $\sigma_\ve\to0$ as $\ve\to0$.
\end{Lem}
\begin{prooff}
Similar to Step.~3 in Lemma~\ref{lem021}, we have
$$
\aligned
I_{\lambda,\mu}(v_{\lambda,\mu}^++t_{\lambda,\mu,\ve}U_{\ve,\overrightarrow{y}})=&I_{\lambda,\mu}(v_{\lambda,\mu}^+)+
I_{\lambda,0}(t_{\lambda,\mu,\ve}U_{\ve,\overrightarrow{y}})
-\frac1{2^*}\int_\Omega|v_{\lambda,\mu}^++\mu\varphi+t_{\lambda,\mu,\ve}U_{\ve,\overrightarrow{y}}|^{2^*}dx\\
&+\frac1{2^*}\int_\Omega|t_{\lambda,\mu,\ve}U_{\ve,\overrightarrow{y}}|^{2^*}dx+\frac1{2^*}\int_\Omega|v_{\lambda,\mu}^++\mu\varphi|^{2^*}dx\\
&+\int_\Omega(v_{\lambda,\mu}^++\mu\varphi)^{2^*-1}t_{\lambda,\mu,\ve}U_{\ve,\overrightarrow{y}}dx.
\endaligned
$$
Since $t_{\lambda,\mu,\ve}\in(t_{\lambda,\mu}', t_{\lambda,\mu}^*)$, by a famous estimate (cf. Sun and Li \cite{SL08}), we have
$$
\aligned
I_{\lambda,\mu}(v_{\lambda,\mu}^++t_{\lambda,\mu,\ve}U_{\ve,\overrightarrow{y}})\leq& m_{\lambda,\mu}^++\frac{t_{\lambda,\mu,\ve}^2}{2}(\|U_{\ve,\overrightarrow{y}}\|^2-\lambda\|U_{\ve,\overrightarrow{y}}\|_2^2)
-\frac{t_{\lambda,\mu,\ve}^{2^*}}{2^*}\|U_{\ve,\overrightarrow{y}}\|_{2^*}^{2^*}\\
&-(t_{\lambda,\mu}')^{2^*-1}\int_\Omega(U_{\ve,\overrightarrow{y}})^{2^*-1}(v_{\lambda,\mu}^++\mu\varphi)dx+o(\ve^{(N-2)/2})\\
\leq&m_{\lambda,\mu}^++\frac1N S^{N/2}-\ve^{(N-2)/2}(D-o(1)),
\endaligned
$$
where
$$
D=(t_{\lambda,\mu}')^{2^*-1}(v_{\lambda,\mu}^+(\overrightarrow{y})+
\mu\varphi(\overrightarrow{y}))\int_{\bbrn}\frac{(N(N-2))^{\frac{N+2}{2}}}{(1+|x|^2)^{(N+2)/2}}dx
\quad\text{for }N\geq4.
$$
Therefore, by Lemma~\ref{lem013}, there exists $\ve^*\leq\ve_0$ such that for $\ve\in(0, \ve^*)$, $v_{\lambda,\mu}^++t_{\lambda,\mu,\ve}U_{\ve,\overrightarrow{y}}\in\mathcal{G}_{\lambda,\mu,\sigma_\ve}$ for all $\overrightarrow{y}\in\mathcal{S}^{N-1}$, where $\sigma_\ve=\ve^{(N-2)/2}(D-o(1))\to0$ as $\ve\to0$.
\end{prooff}

Now, we can obtain the third solution of $(\mathcal{Q}_{\lambda,\mu})$.
\begin{Lem}\label{lem015}
Assume $\lambda\in(0, \lambda_{**}^1)$ and $\mu\in(0, \mu_{**}^1)$.  Then $(\mathcal{Q}_{\lambda,\mu})$ has three solutions.
\end{Lem}
\begin{prooff}
Since Lemmas~\ref{lem012}--\ref{lem014} holds for $\lambda\in(0, \lambda_{**}^1)$ and $\mu\in(0, \mu_{**}^1)$, we can follow the proof of \cite[Lemma~5.6]{W10} step by step to show that Cat$(\mathcal{G}_{\lambda,\mu})\geq2$ by using Lemma~\ref{lem010}.  Thanks to Lemma~\ref{lem009}, there exist two solutions of $(\mathcal{Q}_{\lambda,\mu})$ in $\mathcal{G}_{\lambda,\mu}$.  This, together with Lemmas~\ref{lem003}, \ref{lem006} and \ref{lem007}, implies that $(\mathcal{Q}_{\lambda,\mu})$ has three solutions for $\lambda\in(0, \lambda_{**}^1)$ and $\mu\in(0, \mu_{**}^1)$.
\end{prooff}

\section{The fourth solution of $(\mathcal{Q}_{\lambda,\mu})$}
In this section, we will follow the strategy of He and Yang \cite{HY09} to discuss the existence of the fourth solution for $(\mathcal{Q}_{\lambda,\mu})$.  We begin with
\begin{Lem}\label{lem101}
For every $\lambda\in(0, \lambda_1)$, there exist $\mu_\lambda^{**}\in(0, \mu_\lambda^*)$ and $r_{\lambda}>0$ such that
\begin{itemize}
\item[$(1)$] $I_{\lambda,\mu}(v)$ is strictly convex in $B(0, r_\lambda)$, where $B(0, r):=\{u\in\h: \|u\|<r\}$.
\item[$(2)$] $\mathcal{N}_{\lambda,\mu}^+\subset B(0, r_\lambda)$ for $\mu\in(0, \mu_\lambda^{**})$.
\item[$(3)$] $v_{\lambda,\mu}^+$ is the unique critical point of $I_{\lambda,\mu}$ in $\mathcal{N}_{\lambda,\mu}^+$.
\end{itemize}
\end{Lem}
\begin{prooff} \quad$(1)$\quad By a direct calculation, we have
    $$
    \aligned
    I_{\lambda,\mu}''(u)(v,v)=&\|v\|^2-\lambda\|v\|_2^2-(2^*-1)\int_{\Omega}|u+\mu\varphi|^{2^*-2}v^2dx\\
    \geq&(1-\frac{\lambda}{\lambda_1})S\|v\|_{2^*}^2-(2^*-1)2^{2^*-2}(S^{\frac{2^*-2}{2}}\|u\|^{2^*-2}
    +C\mu^{2^*-2})\|v\|_{2^*}^2\\
    =&\big((1-\frac{\lambda}{\lambda_1})S-(2^*-1)2^{2^*-2}(S^{\frac{2^*-2}{2}}\|u\|^{2^*-2}
    +C\mu^{2^*-2})\big)\|v\|_{2^*}^2.
    \endaligned
    $$
Therefore, for every $\lambda\in(0, \lambda_1)$, there exists $\mu_\lambda^{00}\in(0, \mu_\lambda^*)$ such that
$$
(1-\frac{\lambda}{\lambda_1})S-(2^*-1)2^{2^*-2}C\mu^{2^*-2}>\frac12(1-\frac{\lambda}{\lambda_1})S\quad\text{for }\mu\in(0, \mu_\lambda^{00}).
$$
Take
$$
r_\lambda=\big(\frac{\frac12(1-\frac{\lambda}{\lambda_1})S}{(2^*-1)2^{2^*-2}S^{\frac{2^*-2}{2}}}\big)^{1/(2^*-2)},
$$
then $I_{\lambda,\mu}''(u)(v,v)>0$ for $u\in B(0, r_\lambda)$, $\mu\in(0, \mu_\lambda^{00})$ and $v\in\h\backslash\{0\}$.

$(2)$\quad Similar to \eqref{eq008}, we can show that $\|v\|\leq C\big((\frac{\lambda_1}{\lambda_1-\lambda})((\lambda\mu_n)^{\frac{2^*}{2^*-1}}+\mu_n^{2^*})\big)^{\frac12}$ for all $v\in\mathcal{N}_{\lambda,\mu}^+$.  It follows that there exists $\mu_\lambda^{**}\in(0, \mu_\lambda^{00})$ such that $\mathcal{N}_{\lambda,\mu}^+\subset B(0, r_\lambda)$ for $\mu\in(0, \mu_\lambda^{**})$.

$(3)$\quad Assume a contradiction, then there exist at least two critical points, denoted by $v_{\lambda,\mu}^1$ and $v_{\lambda,\mu}^2$, lie in $\mathcal{N}_{\lambda,\mu}^+$.  By $(1)$, we know that $I_{\lambda,\mu}''(tv_{\lambda,\mu}^1+(1-t)v_{\lambda,\mu}^2)(w,w)>0$ for $\lambda\in(0, \lambda_1)$, $\mu\in(0, \mu_\lambda^{**})$, $w\in\h\backslash\{0\}$ and $t\in[0, 1]$.  This, together with the Taylor's expansion, implies both $I_{\lambda,\mu}(v_{\lambda,\mu}^1)>I_{\lambda,\mu}(v_{\lambda,\mu}^2)$ and $I_{\lambda,\mu}(v_{\lambda,\mu}^2)>I_{\lambda,\mu}(v_{\lambda,\mu}^1)$, which is a contradiction.
\end{prooff}

The following lemma is a further local compactness lemma.
\begin{Lem}\label{lem102}
Assume $v_n\geq0$, $v_n\in\mathcal{N}_{\lambda,\mu}^-$, $I_{\lambda,\mu}(v_n)=c+o_n(1)$ and $I_{\lambda,\mu}'(v_n)=o_n(1)$.  If $c\in(m_{\lambda,\mu}^++\frac1N S^{N/2}, m_{\lambda,\mu}^-+\frac1N S^{N/2})$, then there exists $v_0\in\mathcal{N}_{\lambda,\mu}^-$ such that $v_n\to v_0$ as $n\to\infty$ up to a subsequence for $\lambda\in(0, \lambda_1)$ and $\mu\in(0, \mu_\lambda^{**})$.
\end{Lem}
\begin{prooff}
As in Lemma~\ref{lem004}, we can show that there exists $v_0\in\h$ such that $v_n\rightharpoonup v_0$ in $\h$.  Let $w_n=v_n-v_0$.  Then by \eqref{eq009} and the Br\'ezis--Lieb Lemma, we have
$$
I_{0,0}(w_n)=I_{\lambda,\mu}(v_n)-I_{\lambda,\mu}(v_0)+o_n(1).
$$
Similar to \eqref{eq900}, we can obtain that
$$
I_{0,0}'(w_n)\psi=I_{\lambda,\mu}'(v_n)\psi-I_{\lambda,\mu}'(v_0)\psi+o_n(1)\quad\text{for all }\psi\in\h.
$$
Clearly, $v_0$ is a critical point of $I_{\lambda,\mu}(v)$.  Hence, by a result of Struwe \cite{S00}, there exist $l\in\bbn$ and a solution $v^*$ of $(\mathcal{Q}_{0,0})$ such that $I_{0,0}(w_n)=I_{0,0}(v^*)+\frac lNS^{N/2}+o_n(1)$.  If $v_0\in\mathcal{N}_{\lambda,\mu}^+$, then by Lemmas~\ref{lem003}, \ref{lem021} and \ref{lem101}, we know that $v^*\not=0$ and $l=0$, since it is well-known that $v^*\not=0$ implies $I_{0,0}(v^*)\geq\frac1N S^{N/2}$.  Using Struwe's result again, we known that $v_n\to v_0+v^*$ in $\h$, which is impossible since $v_n\rightharpoonup v_0$ in $\h$.  Therefore, we must have $v_0\in\mathcal{N}_{\lambda,\mu}^-$.  This implies $I_{0,0}(w_n)<\frac 1NS^{N/2}+o_n(1)$.  Since $I_{0,0}(v^*)\geq\frac 1NS^{N/2}$ if $v^*\not=0$, we must have $v^*=0$ and $l=0$.  Due to Struwe's result again, we know that $v_n\to v_0$ in $\h$.
\end{prooff}

In what follows, we re-denote $U_{\ve,\overrightarrow{y}}(x)$ and $\phi(x)$ by $U_{\ve,\overrightarrow{y}}^{\delta_0}(x)$ and $\phi_{\delta_0}(x)$ for the sake of clarity.  Similar to \cite{HY09},  we denote $\mathcal{V}:=\{v\in\h: v\geq0, \|v\|_{2^*}=1\}$.  We define a functional $\beta:\mathcal{V}\to\bbrn$ given by $\beta(v)=\int_{\bbrn}x|v|^{2^*}dx$, where the function $v$ is extend to $\bbrn$ by setting $v=0$ outside $\Omega$.
Set $\mathcal{A}_0:=\{v\in\mathcal{V}: \beta(v)=0\}$ and $c_*=\inf_{v\in\mathcal{A}_0}\|v\|^2$, then by \cite[Lemma~5.4]{HY09}, $c_*>S$.  Thanks to \cite[Lemma~5.6]{HY09}, there exists $\ve^{**}\in(0, \ve^*)$ such that $\|U_{\ve,\overrightarrow{y}}^{\delta_0}(x)\|_{D^{1,2}(\bbrn)}^2\in(S, \frac{S+c_*}{2})$ for $\ve\in(0, \ve^{**})$.  Let $\overline{r}_*=1-\ve^{**}$ and $\mathcal{B}_{\overline{r}_*}:=\{(1-\ve)\overrightarrow{y}\in\bbrn: |(1-\ve)\overrightarrow{y}|\leq\overline{r}_*, \overrightarrow{y}\in\mathcal{S}^{N-1}, 0<\ve<1\}$.  We also define a functional $J_{\lambda,\mu}:\mathcal{V}\to\bbr$ given by $J_{\lambda,\mu}(v)=I_{\lambda,\mu}(t^-_{\lambda,\mu}(v)v)$, where $t^-_{\lambda,\mu}(v)$ is given in Lemma~\ref{lem001}.  Set
$$
\gamma_{\lambda,\mu}=\inf_{\mathcal{F}}\sup_{\mathcal{B}_{\overline{r}_*}}J_{\lambda,\mu}(v),
$$
where $\mathcal{F}:=\{h\in C(\mathcal{B}_{\overline{r}_*}, \mathcal{V}): h|_{\partial \mathcal{B}_{\overline{r}_*}}=\widetilde{U}_{\ve,\overrightarrow{y}}^{\delta_0}(x)\}$ and $\widetilde{U}_{\ve,\overrightarrow{y}}^{\delta_0}=\frac{U_{\ve,\overrightarrow{y}}^{\delta_0}}{\|U_{\ve,\overrightarrow{y}}^{\delta_0}\|_{2^*}}$.  The next three lemmas will help us to show that $\gamma_{\lambda,\mu}$ is a critical value of $J_{\lambda,\mu}$ on $\mathcal{V}$.

\begin{Lem}\label{lem103}
$J_{\lambda,\mu}(\widetilde{U}_{\ve,\overrightarrow{y}}^{\delta_0})=\frac1N S^{N/2}+O(\ve)+O(\mu)$ for every $\delta_0\in(0, 1)$, $\lambda\in(0, \lambda_1)$ and uniformly for $\overrightarrow{y}\in\mathcal{S}^{N-1}$, where
$O(\ve)\to0$ as $\ve\to0$ and $O(\mu)\to0$ as $\mu\to0$.
\end{Lem}
\begin{prooff}
By \cite[Lemma~4.2]{HY09}, $\|U_{\ve,\overrightarrow{y}}^{\delta_0}\|^2=S^{N/2}+o(\ve)$, $\|U_{\ve,\overrightarrow{y}}^{\delta_0}\|_{2^*}^{2^*}=S^{N/2}+o(\ve)$ and $U_{\ve,\overrightarrow{y}}^{\delta_0}\rightharpoonup0$ as $\ve\to0$ in $\h$ for every $\delta_0\in(0, 1)$ and uniformly for $\overrightarrow{y}\in\mathcal{S}^{N-1}$.  These, together with $t^-_{\lambda,\mu}(\widetilde{U}_{\ve,\overrightarrow{y}}^{\delta_0})\widetilde{U}_{\ve,\overrightarrow{y}}^{\delta_0}\in\mathcal{N}_{\lambda,\mu}^-$, implies
$$
\aligned
0=&t^-_{\lambda,\mu}(\widetilde{U}_{\ve,\overrightarrow{y}}^{\delta_0})\|\widetilde{U}_{\ve,\overrightarrow{y}}^{\delta_0}\|^2-\lambda t^-_{\lambda,\mu}(\widetilde{U}_{\ve,\overrightarrow{y}}^{\delta_0})\|\widetilde{U}_{\ve,\overrightarrow{y}}^{\delta_0}\|_2^2-
\lambda\mu\int_\Omega\varphi \widetilde{U}_{\ve,\overrightarrow{y}}^{\delta_0}dx\\
&-\int_\Omega(t^-_{\lambda,\mu}(\widetilde{U}_{\ve,\overrightarrow{y}}^{\delta_0})\widetilde{U}_{\ve,\overrightarrow{y}}^{\delta_0}
+\mu\varphi)^{2^*-1}\widetilde{U}_{\ve,\overrightarrow{y}}^{\delta_0}dx\\
=&t^-_{\lambda,\mu}(\widetilde{U}_{\ve,\overrightarrow{y}}^{\delta_0})(\|\widetilde{U}_{\ve,\overrightarrow{y}}^{\delta_0}\|^2-\lambda
\|\widetilde{U}_{\ve,\overrightarrow{y}}^{\delta_0}\|_2^2-(t^-_{\lambda,\mu}(\widetilde{U}_{\ve,\overrightarrow{y}}^{\delta_0}))^{2^*-2}
)+O(\ve)+O(\mu).
\endaligned
$$
Similar to \eqref{eq007}, for every $\delta_0\in(0, 1)$ and $\overrightarrow{y}\in\mathcal{S}^{N-1}$, $t^-_{\lambda,\mu}(\widetilde{U}_{\ve,\overrightarrow{y}}^{\delta_0})\geq (C+O(\ve))(1-\frac{\lambda}{\lambda_1}+O(\mu))^{\frac{1}{2^*-2}}$.  Hence, we must have
$$
\|\widetilde{U}_{\ve,\overrightarrow{y}}^{\delta_0}\|^2-\lambda
\|\widetilde{U}_{\ve,\overrightarrow{y}}^{\delta_0}\|_2^2-(t^-_{\lambda,\mu}(\widetilde{U}_{\ve,\overrightarrow{y}}^{\delta_0}))^{2^*-2}
=O(\ve)+O(\mu),
$$
which means that $t^-_{\lambda,\mu}(\widetilde{U}_{\ve,\overrightarrow{y}}^{\delta_0})=S^{1/(2^*-2)}+O(\ve)+O(\mu)$.  It follows that
$$
\aligned
J_{\lambda,\mu}(\widetilde{U}_{\ve,\overrightarrow{y}}^{\delta_0})=&\frac12(S^{(N-2)/2}+O(\ve)+O(\mu))(\frac{S^{N/2}+o(\ve)}{S^{(N-2)/2}+o(\ve)})
-\frac{\lambda(O(\ve)+O(\mu))}{2(S^{(N-2)/2}+o(\ve))}\\
&-\frac1{2^*}(S^{N/2}+O(\ve)+O(\mu))\\
=&\frac1N S^{N/2}+O(\ve)+O(\mu)
\endaligned
$$
for every $\delta_0\in(0, 1)$, $\lambda\in(0, \lambda_1)$ and uniformly for $\overrightarrow{y}\in\mathcal{S}^{N-1}$.
\end{prooff}

\begin{Lem}\label{lem104}
$m_{\lambda,\mu}^-=\frac1N S^{N/2}+O(\lambda)+O(\mu)$, where $O(\lambda)\to0$ as $\lambda\to0$ and $O(\mu)\to0$ as $\mu\to0$.
\end{Lem}
\begin{prooff}
By Lemmas~\ref{lem003} and \ref{lem021}, it is easy to see that $m_{\lambda,\mu}^-<\frac1N S^{N/2}$.  Since $v^-_{\lambda,\mu}\in\mathcal{N}_{\lambda,\mu}^-$ is a solution of $(\mathcal{Q}_{\lambda,\mu})$, similar to \eqref{eq008} and \eqref{eq006}, we can conclude that
$$
(C(1-\frac{\lambda}{\lambda_1})-C_1\mu^{2^*-2}))^{1/(2^*-2)}\leq\|v^-_{\lambda,\mu}\|_{2^*}\leq C+C(\lambda,\mu),
$$
where $C(\lambda,\mu)\to0$ as $\lambda\to0$ and $\mu\to0$.  It is well-known that there exists $t^0_{\lambda,\mu}>0$ such that $t^0_{\lambda,\mu}v^-_{\lambda,\mu}\in\mathcal{N}_{0,0}$.  Since $(C(1-\frac{\lambda}{\lambda_1})-C_1\mu^{2^*-2}))^{1/(2^*-2)}\leq\|v^-_{\lambda,\mu}\|_{2^*}$, there exists $T>0$ such that $t^0_{\lambda,\mu}\leq T$ for $\lambda$ and $\mu$ small enough.  It follows from Remark~\ref{rmk001} that
$$
m_{\lambda,\mu}^-=I_{\lambda,\mu}(v^-_{\lambda,\mu})\geq I_{\lambda,\mu}(t^0_{\lambda,\mu}v^-_{\lambda,\mu})\geq\frac1N S^{N/2}+O(\lambda)+O(\mu),
$$
which completes the proof of this lemma.
\end{prooff}

\begin{Lem}\label{lem105}
There exist $\rho>0$, $\delta_*\in(0, 1)$, $\lambda_{**}^2\in(0, \lambda_{**}^1)$ and $\mu_{**}^2\in(0, \mu_{**}^1)$ such that $\frac1N S^{N/2}<\gamma_{0,0}-\rho<\gamma_{\lambda,\mu}<\gamma_{0,0}+\rho<\frac1N S^{N/2}+m_{\lambda,\mu}^-$ for $\lambda\in(0, \lambda_*)$, $\mu\in(0, \mu_*)$ and $\delta_0\in(0, \delta_*)$.
\end{Lem}
\begin{prooff}
Similar to \cite[Lemma~5.10]{HY09}, we can prove that there exists $\delta_*\in(0, 1)$ such that $\gamma_{0,0}\in(\frac1N S^{N/2}, \frac2N S^{N/2})$ for $\delta_0\in(0, \delta_*)$. It follows that there exists $\rho>0$ such that $\frac1N S^{N/2}<\gamma_{0,0}-\rho<\gamma_{0,0}+\rho<\frac2N S^{N/2}$.  Thanks to Lemma~\ref{lem104}, there exist $\lambda^1\in(0, \lambda_{**}^1)$ and $\mu^1\in(0, \mu_{**}^1)$ such that $\gamma_0+\rho<\frac1N S^{N/2}+m_{\lambda,\mu}^-$ for $\lambda\in(0, \lambda^1)$ and $\mu\in(0, \mu^1)$.  On the other hand, for every $d\in(0, 1)$ and $v\in\mathcal{V}$, by the Young inequality, we have
$$
\aligned
J_{\lambda,\mu}(v)=I_{\lambda,\mu}(t^-_{\lambda,\mu}(v)v)\leq& I_{0,0}(t^-_{\lambda,\mu}(v)v)+d\|t^-_{\lambda,\mu}(v)v\|_{2^*}^{2^*}+C(\lambda,\mu,d)\\
\leq&I_{0,0}(t^+_{d}(v)v)+d\|t^+_{d}(v)v\|_{2^*}^{2^*}+C(\lambda,\mu,d)\\
=&(\frac{1}{1-d})^{\frac{N-2}{2}}\|v\|^2+C(\lambda,\mu,d)\\
=&(\frac{1}{1-d})^{\frac{N-2}{2}}J_{0,0}(v)+C(\lambda,\mu,d),
\endaligned
$$
and
$$
\aligned
J_{\lambda,\mu}(v)=I_{\lambda,\mu}(t^-_{\lambda,\mu}(v)v)\geq&
I_{\lambda,\mu}(t^-_{d}(v)v)\\
\geq&I_{0,0}(t^-_{d}(v)v)-d\|t^-_{d}(v)v\|_{2^*}^{2^*}-C(\lambda,\mu,d)\\
=&(\frac{1}{1+d})^{\frac{N-2}{2}}\|v\|^2-C(\lambda,\mu,d)\\
=&(\frac{1}{1+d})^{\frac{N-2}{2}}J_{0,0}(v)-C(\lambda,\mu,d),
\endaligned
$$
where $C(\lambda,\mu,d)\to0$ as $\lambda\to0$ and $\mu\to0$ for fixed $d\in(0, 1)$, $I_{0,0}(t^+_{d}(v)v)+d\|t^+_{d}(v)v\|_{2^*}^{2^*}=\max_{t\geq0}(I_{0,0}(tv)+d\|tv\|_{2^*}^{2^*})$ and $I_{0,0}(t^-_{d}(v)v)-d\|t^-_{d}(v)v\|_{2^*}^{2^*}=\max_{t\geq0}(I_{0,0}(tv)-d\|tv\|_{2^*}^{2^*})$.  By the definition of $\gamma_{\lambda,\mu}$, we have $(\frac{1}{1+d})^{\frac{N-2}{2}}\gamma_{0,0}-C(\lambda,\mu,d)\leq\gamma_{\lambda,\mu}\leq(\frac{1}{1-d})^{\frac{N-2}{2}}\gamma_{0,0}+C(\lambda,\mu,d)$.  It follows that there exist $\lambda_{**}^2\in(0, \lambda^1)$ and $\mu_{**}^2\in(0, \mu^1)$ such that $\gamma_0-\rho<\gamma_{\lambda,\mu}<\gamma_0+\rho$, which completes the proof of this lemma.
\end{prooff}

Combining Lemmas~\ref{lem103} and \ref{lem105}, there exist $\lambda_{*}\in(0, \lambda_{**}^2)$ and $\mu_{*}\in(0, \mu_{**}^2)$ such that $\gamma_{\lambda,\mu}>J_{\lambda,\mu}(\widetilde{U}_{\ve,\overrightarrow{y}}^{\delta_0})$ for $\lambda\in(0, \lambda_{*})$, $\mu\in(0, \mu_{*})$ and $\ve$ small enough.  By the minimax principle (see Ambrosetti and Rabinowitz \cite{AR73}), $\gamma_{\lambda,\mu}$ is a critical value of $J_{\lambda,\mu}(v)$ on $\mathcal{V}$, that is, there exists $\{v_n\}\subset\mathcal{V}$ such that
$$
J_{\lambda,\mu}(v_n)=\gamma_{\lambda,\mu}+o_n(1),\quad \|J'_{\lambda,\mu}(v_n)\|_{T^*_{v_n}\mathcal{V}}=o_n(1),
$$
where $T^*_{v_n}\mathcal{V}$ is the dual space of $T_{v_n}\mathcal{V}:=\{w\in\h: \int_{\Omega}(v_n)^{2^*-1}w=0\}$.
In what follows, we will show the existence of the fourth solution of $(\mathcal{Q}_{\lambda,\mu})$.
\begin{Lem}\label{lem106}
There exists $v_{\lambda,\mu}^4\in\mathcal{N}_{\lambda,\mu}^-$ such that $I_{\lambda,\mu}(v_{\lambda,\mu}^4)=\gamma_{\lambda,\mu}$ and $I'_{\lambda,\mu}(v_{\lambda,\mu}^4)=0$ for $\lambda\in(0, \lambda_*)$ and $\mu\in(0, \mu_*)$.
\end{Lem}
\begin{prooff}
By the implicit function theorem, for $h\in T_{v_n}\mathcal{V}$, we have
$$
J_{\lambda,\mu}'(v_n)h=((t^-_{\lambda,\mu})'(v_n)h)I_{\lambda,\mu}'(t^-_{\lambda,\mu}(v_n)v_n)v_n
+t^-_{\lambda,\mu}(v_n)I_{\lambda,\mu}'(t^-_{\lambda,\mu}(v_n)v_n)h.
$$
By Lemma~\ref{lem001}, for every $v_n$, $t_{\lambda,\mu}^-(v_n)>0$.  This implies $I_{\lambda,\mu}'(t^-_{\lambda,\mu}(v_n)v_n)v_n=0$, since $t^-_{\lambda,\mu}(v_n)v_n\in\mathcal{N}_{\lambda,\mu}^-$.  It follows that
\begin{equation}\label{eq901}
J_{\lambda,\mu}'(v_n)h=t^-_{\lambda,\mu}(v_n)I_{\lambda,\mu}'(t^-_{\lambda,\mu}(v_n)v_n)h\ \ \text{for }h\in T_{v_n}\mathcal{V}.
\end{equation}
We claim that $\{t^-_{\lambda,\mu}(v_n)\}_{n\in\bbn}$ stays positive and bounded away from $0$.  If not, then $t^-_{\lambda,\mu}(v_n)\to0$ up to a subsequence as $n\to\infty$.  It follows that $t^-_{\lambda,\mu}(v_n)(\|v_n\|^2-\lambda\|v_n\|_2^2)=\lambda\mu\int_{\Omega}\varphi v_ndx+\int_{\Omega}(\mu\varphi)^{2^*-1}v_n+o_n(1)$, since $t^-_{\lambda,\mu}(v_n)v_n\in\mathcal{N}_{\lambda,\mu}^-$ and $v_n\in\mathcal{V}$.  Hence,
$$
\frac1N S^{N/2}+o_n(1)\leq\gamma_{\lambda,\mu}+o_n(1)=J_{\lambda,\mu}(v_n)\leq-\frac{\lambda}{2}\|\mu\varphi\|_2^2+o_n(1)<o_n(1),
$$
a contradiction.  Note that $\h=\bbr v_n\oplus T_{v_n}\mathcal{V}$ and $I_{\lambda,\mu}'(t^-_{\lambda,\mu}(v_n)v_n)v_n=0$,
$$
I_{\lambda,\mu}'(t^-_{\lambda,\mu}(v_n)v_n)w=I_{\lambda,\mu}'(t^-_{\lambda,\mu}(v_n)v_n)h_w\quad\text{for all }w\in\h,
$$
where $h_w$ is the projection of $w$ in $T_{v_n}\mathcal{V}$.  Thus, by \eqref{eq901},
$$
\aligned
\|I_{\lambda,\mu}'(t^-_{\lambda,\mu}(v_n)v_n)\|=&\sup_{w\in\h\backslash\{0\}}\frac{I_{\lambda,\mu}'(t^-_{\lambda,\mu}(v_n)v_n)w}{\|w\|}\\
\leq&\sup_{w\in\h\backslash\{0\}}\frac{I_{\lambda,\mu}'(t^-_{\lambda,\mu}(v_n)v_n)h_w}{\|h_w\|}\\
=&\sup_{w\in\h\backslash\{0\}}\frac{J_{\lambda,\mu}'(v_n)h_w}{t^-_{\lambda,\mu}(v_n)h_w}\\
\leq& C\|J'_{\lambda,\mu}(v_n)\|_{T^*_{v_n}\mathcal{V}}=o_n(1).
\endaligned
$$
On the other hand, by the definition of $J_{\lambda,\mu}$, it is easy to see that $I_{\lambda,\mu}(t^-_{\lambda,\mu}(v_n)v_n)\to\gamma_{\lambda,\mu}$ as $n\to\infty$.  Thanks to Lemmas~\ref{lem102} and ~\ref{lem105}, $t^-_{\lambda,\mu}(v_n)v_n\to v_{\lambda,\mu}^4$ as $n\to\infty$ for $\lambda\in(0, \lambda_*)$, $\mu\in(0, \mu_*)$ and $\delta_0\in(0, \delta_*)$.  Clearly, $v_{\lambda,\mu}^4$ is a solution of $(\mathcal{Q}_{\lambda,\mu})$ and different from the solutions find in Lemmas~\ref{lem006}, \ref{lem007} and \ref{lem015} for $\lambda\in(0, \lambda_*)$, $\mu\in(0, \mu_*)$ and $\delta_0\in(0, \delta_*)$.
\end{prooff}

We close this section by
\medskip\par\noindent{\bf Proof of Theorem~\ref{thm004}}\quad It follows immediately from Lemmas~\ref{lem015} and \ref{lem106}.
\qquad\raisebox{-0.5mm}{\rule{1.5mm}{4mm}}\vspace{6pt}

\section{Acknowledgements}
Y. Wu is supported by the Fundamental Research Funds for the Central Universities (2014QNA67).

\end{document}